\documentclass[11pt]{article}
\usepackage{amsmath,amssymb}
\usepackage{times}
\usepackage{bm}
\usepackage{algorithm}
\usepackage{amssymb}
\usepackage{mathrsfs}
\usepackage{graphicx}
\usepackage[flushleft]{threeparttable}
\usepackage{enumitem} 
\usepackage[normalem]{ulem}

\usepackage{star}
\usepackage[colorlinks,
linkcolor=blue,
anchorcolor=blue,
citecolor=blue
]{hyperref}

\usepackage{array}
\usepackage{makecell}

\newtheorem{theorem}{Theorem}[section]
\DeclareRobustCommand\optionalsec[1]{%
  \ifnum\pdfstrcmp{#1}{\thesection}=0\else#1.\fi
}

\newtheorem{lemma}{Lemma}[section]
\DeclareRobustCommand\optionalsec[1]{%
  \ifnum\pdfstrcmp{#1}{\thesection}=0\else#1.\fi
}

\theoremstyle{remark}
\newtheorem{remark}{Remark}[subsection]

\newcommand {\Sym}{{\rm Sym}}

\def\##1\#{\begin{align}#1\end{align}}
\def\$#1\${\begin{align*}#1\end{align*}}

\newcommand{\Rom}[1]{\text{\uppercase\expandafter{\romannumeral #1\relax}}}
\usepackage{txfonts}
\usepackage{geometry}
\begin{document}

\title{Restricted Riemannian geometry for positive semidefinite matrices}  

\author{A. Martina Neuman\thanks{Department of Computational Mathematics Science and Engineering, Michigan State University, 428 S Shaw Lane, East Lansing, 48824, USA; E-mail: \texttt{neumana6@msu.edu}.}, Yuying Xie \thanks{Department of Computational Mathematics Science and Engineering and Department of Statistics and Probability, Michigan State University, 428 S Shaw Lane, East Lansing, 48824, USA; E-mail: \texttt{xyy@msu.edu}.}, and Qiang Sun\thanks{Department of Statistical Sciences, University of Toronto. 100 St. George Street, Toronto, ON M5S 3G3, Canada; E-mail: \texttt{qsun@utstat.toronto.edu} } 
}

\date{ }

\maketitle

\vspace{-0.25in}

% typeset the title of the contribution
\begin{abstract}
We introduce the manifold of {\it restricted} $n\times n$ positive semidefinite matrices of fixed rank $p$, denoted $S(n,p)^{*}$. The manifold itself is an open and dense submanifold of $S(n,p)$, the manifold of $n\times n$ positive semidefinite matrices of the same rank $p$, when both are viewed as manifolds in $\mathbb{R}^{n\times n}$. This density is the key fact that makes the consideration of $S(n,p)^{*}$ statistically meaningful. We furnish $S(n,p)^{*}$ with a convenient, and geodesically complete, Riemannian geometry, as well as a Lie group structure, that permits analytical closed forms for endpoint geodesics, parallel transports, Fr\'echet means, exponential and logarithmic maps. This task is done partly through utilizing a {\it reduced} Cholesky decomposition, whose algorithm is also provided. We produce a second algorithm from this framework to estimate principal eigenspaces and demonstrate its superior performance over other existing algorithms.
\end{abstract}



%\tableofcontents
\section{Introduction}\label{section 1}

\noindent
Positive semidefinite matrices (PSD) of fixed rank appear in many corners of applied mathematics and in many forms, ranging from low-rank covariance matrices (covariance
matrices in statistics \cite{huber2004robust}, sparse ``shape" covariance matrices in action recognition in video \cite{Guo:2010}, covariance descriptors in image set classification \cite{faraki2016image}) to kernel matrices in machine learning \cite{lanckriet2004learning}. Their collective geometry also provides natural frameworks for numerous application tasks. For instance, optimization on these matrices is used in distance matrix completion \cite{mishra2011low} and in computing of the maximal cut of a graph and sparse principal
component analysis \cite{journee2010low}, while in \cite{bakker2018dynamic} and in \cite{li2014conformational}, their geodesic interpolations are utilized in the diverse contexts of protein configuration and dynamic community detection, respectively. This suggests that the investigation into these low-rank structures continues to play an important role in the future.\\

\noindent
Let $n\geq p$ be two positive integers. The set of positive semidefinite $n\times n$ matrices of fixed rank $p$ is denoted by $S(n,p)$. When $n=p$, $S(n,p)=S(n,n)$ is also denoted as $P(n)$ in the literature, and $P(n)$ possesses a well-known conic geometry with a natural non-Euclidean metric \cite{bonnabel2010riemannian}, \cite{faraut1994analysis}. However, when $n>p$, there is no natural choice for geometry for $S(n,p)$ \cite{massart2020quotient}. Instead, a body of work has been dedicated to proposing multiple geometries for $S(n,p)$, each with different applications and goals in mind. Many view $S(n,p)$ as a quotient geometry modulo orthogonal actions, such as $S(n,p)\cong\mathbb{R}^{n\times p}_{*}/O(p)$ in \cite{absil2009geometric}, \cite{journee2010low}, \cite{massart2020quotient}, or $S(n,p)\cong (V(n,p)\times P(p))/O(p)$ in \cite{bonnabel2010riemannian}, \cite{meyer2009subspace} - here, $\mathbb{R}^{n\times p}_{*}, V(n,p)$ are noncompact and compact Stiefel manifolds, respectively \cite{absil2009optimization}. Some others regard $S(n,p)$ as a submanifold embedded in $\mathbb{R}^{n\times n}$ \cite{helmke1995critical}, \cite{koch2007dynamical}, \cite{orsi2006newton}, \cite{vandereycken2009embedded}, etc, while still some like \cite{vandereycken2013riemannian} think of $S(n,p)$ as a homogeneous space of $GL_{+}(n)$, the connected component at the identity of the general linear group. Each of these geometries paints a vivid picture of $S(n,p)$, but none gives a finished one. For instance, \cite{massart2020quotient} provided a notion of distance that coincides with the Wasserstein distance between degenerate, centered Gaussian distributions, but the resulting geometry does not make $S(n,p)$ a complete metric space. \cite{vandereycken2013riemannian} managed to give $S(n,p)$ a geodesically complete geometry with closed form, but not endpoint, geodesics, and expressions for the exponential map, but not for the logarithm. \cite{bonnabel2010riemannian} gives intuitive, physically motivated, notions of distance and geometric-like midpoint, but their distances are not geodesic distances.\\

\noindent
We start with a different viewpoint: we consider a special subset of $S(n,p)$ that is sufficient for statistical applications - a few examples of which are given in \ref{subsubsection 3.4.1} as a motivation. That subset should also be dense in $S(n,p)$ and ``special" enough that it can be endowed with a complete and convenient Riemannian geometry. In a way, it's akin to opting for invertible matrices out of square matrices in developing computing algorithms, as the set of non-invertible matrices has Lebesgue measure zero \cite{strichartz2000way}. Relaxations are not new in applied geometry; for instance, the concept of vector transport was introduced to ease the computational requirements of parallel translation, just as the concept of retraction does the same for exponential mapping \cite{absil2009optimization}. In addition, even though completeness is desirable in many optimization algorithms, there are retraction-based Riemannian optimizations (\cite{adler2002newton}) that do away with the need for geodesics. In the same spirit of relaxation, we consider the subset $S(n,p)^{*}$, called {\it restricted} $S(n,p)$, formally introduced in subsection \ref{subsection 2.1}, that can be given a geodesically complete Riemannian geometry, with analytical forms for endpoint geodesics, exponential and logarithm maps, and parallel transports. Moreover, this set is an open and dense submanifold of $S(n,p)$ with full dimension $np-p(p-1)/2$ \cite{bonnabel2010riemannian}, \cite{massart2020quotient}.\\

\noindent
The theoretical approach of this paper is inspired by that of \cite{lin2019riemannian}, in which the author constructed a Riemannian geometry for the manifold of full-rank symmetric positive definite matrices (denoted $\mathcal{S}^{+}_{m}$ in that paper with $m$ being the matrix size) with all the desirable properties, by constructing the same for the so-called Cholesky space ($\mathcal{L}_{+}$) and mapping it back to $\mathcal{S}^{+}_{m}$ via a manifold diffeomorphism that is the Cholesky decomposition. This technique fails for the low-rank situation, for the Cholesky map is not continuous on $S(n,p)$, when $n>p$, as explained in subsection \ref{subsection 9.2}. However, this map will be continuous on different subsets of $S(n,p)$ if one compromises in advance on the locations of the first $p$ linearly independent columns. The reader is referred to the appendix for more details. In particular, this map turns out to be a diffeomorphism between $S(n,p)^{*}$ and a space similar to the Cholesky space. This approach is also motivated by our observation that a randomly sampled covariance matrix from $S(n,p)$ is always in $S(n,p)^{*}$ (see subsection \ref{subsection 9.3}).\\

\noindent
This Cholesky approach allows one to fall back on the Euclidean nature of the Cholesky space, which is of great asset in deriving simple closed forms for endpoints geodesics, exponential and logarithmic maps, and parallel transports, which in turns aid computation. For instance, since parallel transport is a ubiquitous tool in statistical analysis, machine learning on manifolds \cite{schiratti2017bayesian}, \cite{zeestraten2017approach} or optimizations that involve manifold-valued data \cite{hosseini2015matrix}, it's advantageous in having a closed form to handle big data.\\

\noindent
For applications, we use our derived notion of Fr\'echet mean to estimate the top principal eigenspaces in the setting of distributed PCA.

\subsection{Lay-out}

\noindent
The structure of this paper is as follows. In section \ref{section 2}, as well as the appendix, section \ref{section 9}, we introduce $S(n,p)^{*}$ and the Cholesky-like space $\mathcal{L}^{n\times p}_{*}$ and state our main theorem. In section \ref{section 9}, various needed claims are proved. In section \ref{section 3}, we derive Riemannian geometries for both manifolds and give a closed form for endpoint geodesics as well as for Fr\'echet means on $S(n,p)^{*}$. We give several discussions on applications and illustrations of this geometry at the end of section \ref{section 3}. In section \ref{section 4}, we give expressions for the exponential and logarithm maps, and in section \ref{section 5}, parallel transport. We introduce a quick algorithm to calculate {\it reduced} Cholesky factors and our LRC algorithm for Fr\'echet mean in section \ref{section 6}. We give numerical applications in section \ref{section 7} and finally some concluding remarks in section \ref{section 8}.

\subsection{Notations} \label{subsection 1.2}

\noindent
Let $n>p$ and $k$ be positive integers.

\begin{itemize}
    \item $\mathbb{R}_{*}^{n\times p}$: The manifold of $n\times p$ matrices of full rank, or the noncompact Stiefel manifold of $n\times p$ matrices \cite{absil2009optimization}
    \item $V(n,p)$: The manifold of $n\times p$ matrices of full rank such that $V^{T}V=I_{p}$, or the compact Stiefel manifold of $n\times p$ matrices \cite{absil2009optimization}
    \item $S(n,p)$: The manifold of $n\times n$ positive semidefinite matrices of rank $p$ and of dimension $np-p(p-1)/2$
    \item $\mathcal{L}$: The Cholesky set of $S(n,p)$ (see subsection \ref{subsection 9.2})
    \item $P(k)$: The manifold of $k\times k$ symmetric positive definite matrices with positive eigenvalues \cite{lin2019riemannian}
    \item $O(k)$: The Lie group of $k\times k$ orthogonal matrices \cite{lee2009manifolds}
    \item $\mathcal{P}(k)$: The finite group of $k\times k$ permutation matrices ($\subset O(k)$)
    \item $\Sym(k)$: The linear manifold of $k\times k$ symmetric matrices
    \item $GL_{+}(k)$: The connected component at the identity of the general linear group $GL(k)$ of $k\times k$ invertible matrices \cite{vandereycken2013riemannian}
    \item $\mathcal{D}_{+}(p)$: The set of $p\times p$ diagonal matrices with strictly positive diagonal entries listed decreasingly
    \item $Diag(p)$: The linear manifold of $p\times p$ diagonal matrices
    \item $S(n,p)^{*}$: The manifold of {\it restricted} $S(n,p)$ (see section \ref{section 2})
    \item $\mathcal{L}^{n\times p}_{*}$: The {\it reduced} Cholesky space of $S(n,p)^{*}$ (see subsection \ref{subsection 2.1}, \ref{subsection 9.3})
    \item $\mathcal{L}^{n\times p}$: The vector space of $n\times p$ mock lower triangular matrices (see subsection \ref{subsection 2.1} for the terminology)
\end{itemize}

\section{Preliminaries \& main theorem} \label{section 2}

\noindent
Let $M\in S(n,p)$. As mentioned in section \ref{section 1}, there have been attempts to look at $M$ through a geometric lens. One can write $M=ZZ^{T}$ with $Z\in\mathbb{R}^{n\times p}_{*}$ as in \cite{absil2009geometric}, or $M= UR^2U^{T}$ with $U\in V(n,p), R\in P(p)$ as in \cite{bonnabel2010riemannian}. The problem with these geometric representations is that they are unique up to modulo an orthogonal action. For instance, if $Z\mapsto ZO$, for any $O\in O(p)$, then $(ZO)(ZO)^{T}=ZZ^{T}$, and if $(U,R^2)\mapsto (UO, O^{T}R^2 O)$, then $(UO)(O^{T}R^2 O)(UO)^{T}=UR^2U^{T}$. Even when one restricts attention to $M\in S(n,p)$ whose positive spectrum (positive eigenvalues) is distinct and na\"ively hopes to write 
\begin{equation} \label{diagdec} 
    M=VD^2V^{T},
\end{equation} 
with $D\in\mathcal{D}_{+}(p)$ and $V\in V(n,p)$, to avoid the problem encountered in \cite{bonnabel2010riemannian} then one still has a problem with non-uniqueness of $V$; for instance $M=(V\mathcal{I})D^2(V\mathcal{I})^{T}$, where $\mathcal{I}\in Diag(p)$ has $\pm 1$ on the diagonal. In addition, endpoint geodesics on $V(n,p)$ are not known \cite{bryner2017endpoint}, which means that even \ref{diagdec} doesn't help with establishing endpoint geodesics on $S(n,p)$. We decided to look at the algebraic structures of $M\in S(n,p)$ instead.\\

\noindent
Roughly speaking, one can classify $M$ according to the locations of its first $p$ linearly independent columns. Let $i_{j}=1,\cdots,n$ where $j=1,\cdots,p$ and $i_{j_1}\not= i_{j_2}$ if $j_1\not=j_2$. Let $i_1,\cdots,i_{p}$ then indicate the indices of the first $p$ linearly independent columns of $M$. The set of all these $\{i_1,\cdots,i_{p}\}$ has cardinality ${n\choose p}$. Let $C_{\{i_1,\cdots,i_{p}\}}$ denote the collection of all matrices in $S(n,p)$ whose first $p$ linearly independent columns are precisely the $i_1th,\cdots,i_{p}th$ columns. Then $M\in S(n,p)$ belongs to one and only one of these disjoint collections $C_{\{i_1,\cdots,i_{p}\}}$.\\ 

\noindent
One special case is when $C_{\{i_1,\cdots,i_{p}\}}=C_{\{1,\cdots,p\}}$; ie, when the first $p$ linearly independent columns of $M$ are also its first $p$ columns. Any matrix in $S(n,p)$ can be ``looped" back to a matrix in $C_{\{1,\cdots,p\}}$. For instance, suppose $n=3, p=2$ and, $M = \begin{pmatrix}1&1&1\\ 1&1&1\\ 1&1&2\end{pmatrix}$; i.e., $M\in C_{\{1,3\}}$. Then there exists $P\in\mathcal{P}(3)$ such that $PMP^{T}=PMP\in S(3,2)^{*}$. Take $P = \begin{pmatrix}1&0&0\\ 0&0&1\\ 0&1&0\end{pmatrix}$ for example, then $PMP= \begin{pmatrix}1&1&1\\ 1&2&1\\ 1&1&1\end{pmatrix}\in C_{\{1,2\}}$. In fact, this ``looping" action was at the heart of the incompatibility between motions in $Gr(n,p)$ and $P(p)$ in \cite{bonnabel2010riemannian}. (In that paper, the authors constructed the quotient geometry $(V(n,p)\times P(p))/O(p)$ for $S(n,p)$ to take advantage of the available endpoint geodesics on $Gr(n,p)$ and $P(p)$; however, their choice of horizontal and ``vertical" spaces turned this quotient into an unorthodox skew projection rather than an orthogonal one. This makes for a situation where $A\in S(n,p)$ can get to $B\in S(n,p)$ much ``faster" through motions in $Gr(n,p)$ instead of in $P(p)$.) However, there are also situations where it doesn't require any action to bring matrices in $S(n,p)$ back to $C_{\{1,\cdots,p\}}$. For instance, any of the following matrices is already in $C_{\{1,2\}}$:
\begin{equation*}
    \begin{pmatrix}1&2&3\\2&5&8\\3&8&13\end{pmatrix},\begin{pmatrix}13&8&3\\8&5&2\\3&2&1\end{pmatrix}, \begin{pmatrix}5&2&8\\2&1&3\\8&3&13\end{pmatrix},\begin{pmatrix}13&3&8\\3&1&2\\8&2&5\end{pmatrix},\begin{pmatrix}1&3&2\\3&13&8\\2&8&5\end{pmatrix},\begin{pmatrix}5&8&2\\8&13&3\\2&3&1\end{pmatrix}.
\end{equation*}
Note that these matrices are all possible permutations of the first matrix and that their all being in $C_{\{1,2\}}$ is partly due to $n-1=p$ in this case. This example thwarts any attempt to write $S(n,p)$ as folded copies of $C_{\{1,\cdots,p\}}$. However, there are still advantages to considering just the set $C_{\{1,\cdots,p\}}$ alone, as we further explain.\\

\noindent
Fix $i_1,\cdots,i_{p}$. It's known that every positive semidefinite matrix has a unique Cholesky factorization \cite{gentle2012numerical}. If $M\in C_{\{i_1,\cdots,i_{p}\}}$ then its Cholesky factor, $L$, is a lower triangular matrix with precisely $p$ positive diagonal entries at $L_{i_1\,i_1},\cdots,L_{i_{p}\,i_{p}}$ and $n-p$ zero columns. The converse is also true; ie, the Cholesky factor $L$ of $M$ has a zero column whenever the corresponding column of $M$ is dependent on the previous columns. The reader is referred to sub-subsection \ref{subsubsection 9.1.1} for a proof of this claim. That means that the shape of the lower triangular matrix $L$ is determined whenever $C_{\{i_1,\cdots,i_{p}\}}\ni M$ is determined. That means, if $\{i_1,\cdots,i_{p}\}=\{1,\cdots,p\}$ and $M\in C_{\{1,\cdots,p\}}$, then the Cholesky factor $L$ of $M$ will look like,
\begin{equation*} 
    L=\begin{pmatrix}a_1&0&\dots&0&0&\dots&0\\ \star&a_2&\dots&0&0&\dots&0\\ \vdots&\vdots&\ddots&\vdots&\vdots&\dots&\vdots\\ \star&\star&\dots&a_{p}&0&\dots&0\\ \vdots&\vdots&\dots&\vdots&\vdots&\dots&\vdots\\ \star&\star&\dots&\star&0&\dots&0\end{pmatrix},
\end{equation*}
where $a_1,\cdots,a_{p}>0$, the last $n-p$ columns are all zero columns, and the entries $\star$ below $a_{i}$'s are real entries. It's easy to see that these Cholesky factor matrices form an open submanifold of $\mathbb{R}^{np-p(p-1)/2}$, with its manifold topology inherited from the Euclidean space. That means, collectively, these matrices $L$ form a manifold of dimension $np-p(p-1)/2$, which is also the dimension of $S(n,p)$ \cite{bonnabel2010riemannian}, \cite{vandereycken2013riemannian}. Since the Cholesky map $M\mapsto L$ is one-to-one (see \cite{gentle2012numerical} and the appendix), this implies that the corresponding matrices $M\in C_{\{1,\cdots,p\}}$ should also form a manifold of dimension $np-p(p-1)/2$. This is the idea that we pursue below.\\

\noindent
It's easily deduced from the discussion above that if $\{i_1,\cdots,i_{p}\}\not=\{1,\cdots,p\}$, then the Cholesky matrices $L$ of $M\in C_{\{i_1,\cdots,i_{p}\}}$ form a manifold of dimension strictly less than $np-p(p-1)/2$. For instance, if $n=4, p=2$ and $M\in C_{\{1,3\}}$, then
\begin{equation*} 
    L=\begin{pmatrix}a_1&0&0&0\\ \star&0&0&0\\ \star&0&a_3&0\\ \star&0&\star&0\end{pmatrix}.
\end{equation*}
Here, $a_1,a_3>0$ and $\star\in\mathbb{R}$. Clearly, these $L$'s form a manifold of dimension $4+2=6< 7=4\times 2-2(2-1)/2$. An implication here is that if $M\in S(n,p)$ then it's most likely that $M\in C_{\{1,\cdots,p\}}$.\\

\noindent
In addition, any $M\in S(n,p)$ can be easily approximated by some $M_{\epsilon}\in C_{\{1,\cdots,p\}}$. For instance, suppose $n=3, p=2$ and, $M=\begin{pmatrix}1&1&1\\1&1&1\\1&1&2\end{pmatrix}$. Then one can take, for some $\epsilon>0$: $M_{\epsilon} = \begin{pmatrix}1&1&1\\ 1&1+\epsilon^2&1+\epsilon\\ 1&1+\epsilon&2\end{pmatrix}\in C_{\{1,2\}}$. Then $\|M_{\epsilon}-M\|_{F}=O(\epsilon)$. This reinforces the denseness observation above.\\

\noindent
All these claims (Cholesky factorization, manifold structures, dimensions, density, approximation) will be further explored and proved throughout the paper, but for now, that gives us enough motivation to consider these special matrices $M\in C_{\{1,\cdots,p\}}\subset S(n,p)$.

\subsection{ $S(n,p)^{*}$ and $\mathcal{L}^{n\times p}_{*}$} \label{subsection 2.1}

\noindent
As indicated above, we will restrict our attention to only the subset $C_{\{1,\cdots,p\}}$ of $S(n,p)$. We rename it to {\it restricted} $S(n,p)$ and denote it $S(n,p)^{*}$. We show in subsection \ref{subsection 9.3} that, instead of considering the full Cholesky factor $L\in\mathbb{R}^{n\times n}$ of $M\in S(n,p)^{*}$, one can consider its {\it reduced} Cholesky factor $N\in\mathbb{R}^{n\times p}_{*}$ instead:
\begin{equation*}
    LL^{T}=M=NN^{T}.
\end{equation*}
The set of all these reduced Cholesky factor matrices constitutes a set that is Cholesky-like in nature; we call it the {\it reduced} Cholesky space and denote it $\mathcal{L}^{n\times p}_{*}$. 

\vspace{5pt}

\begin{remark}
Now $N$ here is a ``mock" lower triangular matrix, meaning, $N_{ij}=0$ if $i<j, i=1,\cdots,p$. In this paper, we will consistently call such non-square lower (upper) triangular matrices - that resemble their square counterparts - mock lower (upper) triangular matrices, and their ``diagonals" - $N_{ii}$ - mock diagonals. We ignore citing the ``upper" half entries $N_{ij}, i<j$, in this paper, as they are all zeros.
\end{remark}

\vspace{5pt}

\noindent
We are now ready to state our theorem. 

\subsection{Main theorem}

\begin{theorem}
Let $n>p$ be two positive integers. Let $S(n,p)^{*}$ be as above. Then, \\
1. $S(n,p)^{*}$ is an open and dense submanifold of $S(n,p)$ with the same dimension $np-p(p-1)/2$, when both are viewed as submanifolds of $\mathbb{R}^{n\times n}$. \\
2. $S(n,p)^{*}$ is also a Riemannian manifold with the Riemannian metric $g$ defined in \ref{pull} and endpoint geodesics \ref{geocurv}. A description of the tangent space $T_{M}S(n,p)^{*}=T_{M}S(n,p)$ is given.\\
3. The Fr\'echet mean of $k$ number of points on $S(n,p)^{*}$ exists and has the form \ref{mean1}. If $k=2$ then this Fr\'echet mean is also their geodesic midpoint. This Fr\'echet mean satisfies permutation invariance.\\
4. $S(n,p)^{*}$ has a Lie group structure, with group operation $\otimes$ defined in \ref{otimes}. The Riemannian metric $g$ is also a bi-variant metric with this group structure.\\
5. $S(n,p)^{*}$ has a simple form of parallel transport in \ref{par}. It also has simple forms of exponential and logarithmic maps, defined in \ref{exp}, \ref{log}, respectively.
\end{theorem}

\noindent
Point 1 is proved in sections \ref{section 9} and \ref{section 3}; points 2, 3 are proved in section \ref{section 3}, and the last two points are in sections \ref{section 4}, \ref{section 5}.

\subsection{General remarks}

\noindent
Some arguments presented here are illustrated with concrete choices of $n,p$ (say $n=3,p=2$). We emphasize here that they are used solely to get the point across and to avoid the cumbersome inductive Cholesky factorization formulas for general $n,p$. We opt for these concrete choices for the sake of presentation and readability. Arguments will be presented to showcase that the logic stays the same for general $n,p$.\\

\noindent
Since the paper is heavy in notations, some decisions are made to promote readability that might inadvertently lead to ambiguity. Ambiguity is sometimes also the result of coincidence in notations from different principles. The notation $\mathcal{F}$ in subsection \ref{subsection 3.3} simply means ``mean" (average); it changes its meaning to mean the Fr\'echet mean in $\mathcal{L}^{n\times p}_{*}$ or in $S(n,p)^{*}$ depending on context. The notation $O(k)$ means the Lie group of orthogonal matrices $\mathbb{R}^{k\times k}$, and $O(\epsilon)$ indicates an upper bound that is at most $C\epsilon$ in magnitude, for some $C>0$ and small $\epsilon$.

\section{The manifolds $\mathcal{L}^{n\times p}_{*}$ and $S(n,p)^{*}$} \label{section 3}

\noindent
We prove in subsection \ref{subsection 9.2} that $S(n,p)^{*}$ is a manifold with dimension $np-p(p-1)/2$ and in subsection \ref{subsection 9.3} that $S(n,p)^{*}$ is an open submanifold of $S(n,p)$. This is done so as our main goal is to establish a Riemannian structure of $S(n,p)^{*}$. We do so in this section. We start out by following the arguments in subsection \ref{subsection 9.2} again: we establish a manifold structure on $\mathcal{L}^{n\times p}_{*}$ and push it onto $S(n,p)^{*}$ via \ref{1to1}, using $N$ instead of $L$. We use a non-Euclidean coordinate chart here.\\

\noindent
It's shown in subsection \ref{subsection 9.3} that if $M\in S(n,p)^{*}$ then there exists a unique $N\in\mathcal{L}^{n\times p}_{*}$, called the reduced Cholesky factor of $M$ such that:
\begin{equation*}
    NN^{T}=M,
\end{equation*}
and that the mock diagonal entries $N_{ii}>0, i=1,\cdots,p$. Let $\Psi: \mathcal{L}^{n\times p}_{*}\to S(n,p)^{*}$ be the map:
\begin{equation} \label{psi1} 
    \Psi(N) = NN^{T}
\end{equation}
and $\Psi^{-1}: S(n,p)^{*}\to\mathcal{L}^{n\times p}_{*}$ the map:
\begin{equation} \label{psi2} 
    \Psi^{-1}(M) = N
\end{equation}
where $N$ satisfies \ref{psi1}; see an example of $\Psi$ in \ref{sampM}. Then $\Psi,\Psi^{-1}$ pair points in $\mathcal{L}_{*}^{n\times p}$ and points in $S(n,p)^{*}$ in a one-to-one manner. Structurally speaking, $\mathcal{L}^{n\times p}_{*}$ lives in a Euclidean space $\mathbb{R}^{n\times p}$, so it has the natural ``flat" Euclidean coordinate chart: $x^{ij}(N)=N_{ij}$, for $1\leq j\leq p, j\leq i$. We follow \cite{lin2019riemannian} by ``slowing" the mock positive diagonal down by putting a log growth on it; precisely, we use, 
\begin{equation} \label{Lcoor} 
    x^{ii}(N) = \log N_{ii}, \,i=1,\cdots,p\,\text{ and }\, x^{ij}(N)=N_{ij},\,i>j.
\end{equation} 
That means, in terms of \ref{psi1}, \ref{psi2} (see also \ref{map1}, \ref{chp}), we can define a corresponding chart for $S(n,p)^{*}$:
\begin{equation} \label{coor}
    \begin{split}
        x^{ii}(M) & = \log N_{ii} = \log\bigg(\sqrt{M_{ii}-\sum_{k=1}^{i-1}N_{ik}^2}\bigg),\,\, i=1,\cdots,p\\
        x^{ij}(M) & = N_{ij} =\frac{1}{N_{ii}}\bigg(M_{ij}-\sum_{k=1}^{j-1}N_{ik}N_{jk}\bigg),\,\, i>j.
    \end{split}
\end{equation}
The reader is noted that, given fixed $n,p$, $N_{ij}$'s are solved completely in terms of $M_{ij}$'s (and so are the right hand sides of \ref{coor}) through summations, multiplications, divisions and taking roots. There are $np-p(p-1)/2$ for all the coordinates $(ij), j=1,\cdots,p, j\leq i$. It's clear that $x$ is a local chart on $S(n,p)^{*}$ and the fact that $S(n,p)^{*}$ is a manifold of dimension $np-p(p-1)/2$ doesn't change.

\subsection{The Riemannian manifold $\mathcal{L}^{n\times p}_{*}$}

\noindent
Recall from subsection \ref{subsection 1.2} that $\mathcal{L}^{n\times p}$ is the vector space of $n\times p$ mock lower triangular matrices (see also subsection \ref{subsection 9.3}). Note that $\mathcal{L}^{n\times p}_{*}$ is an open subset of $\mathcal{L}^{n\times p}$, and both are Euclidean manifolds of the same dimension $np-p(p-1)/2$, and that if $N\in\mathcal{L}^{n\times p}_{*}$ then $T_{N}\mathcal{L}^{n\times p}_{*}\cong\mathcal{L}^{n\times p}$. We use the coordinate chart $x$ defined in \ref{Lcoor} for $\mathcal{L}^{n\times p}_{*}$ and define a Riemannian structure on $\mathcal{L}^{n\times p}_{*}$. Let $X,Y\in T_{N}\mathcal{L}^{n\times p}_{*}\cong\mathcal{L}^{n\times p}$, we define $\tilde{g}$ on $\mathcal{L}^{n\times p}_{*}$ as:
\begin{equation} \label{Rie} 
    \tilde{g}_{N}(X,Y)=\sum_{i>j} X_{ij}Y_{ij} + \sum_{i=1}^{p} N_{ii}^{-2}X_{ii}Y_{ii}.
\end{equation}
Define a curve on $\mathcal{L}^{n\times p}_{*}$, starting with $N\in\mathcal{L}^{n\times p}_{*}$ and ``velocity" $X\in\mathcal{L}^{n\times p}$ as,
\begin{equation} \label{geo} 
    \begin{split}
        \tilde{\gamma}_{N,X}(t)_{ii} &= N_{ii}\exp(tX_{ii}N_{ii}^{-1}),\,\, i=1,\cdots,p\\
        \tilde{\gamma}_{N,X}(t)_{ij} &= N_{ij}+tX_{ij},\,\, i>j.
    \end{split}  
\end{equation}
It's clear that $\tilde{\gamma}_{N,X}(t)\in\mathcal{L}^{n\times p}_{*}$ for all $t\in(-\infty,\infty)$. As in \cite{lin2019riemannian}, it's easy to show that $\tilde{\gamma}_{N,X}(t)$ is a geodesic for $(\mathcal{L}^{n\times p}_{*},\tilde{g})$ - which then follows that $\mathcal{L}^{n\times p}_{*}$ is geodesically complete. The idea is to show that $\tilde{\gamma}_{N,X}$ satisfies the geodesic equations \cite{lee2009manifolds}, but this follows exactly as in \cite{lin2019riemannian}, so we refer the reader to the said paper. We state here a point about the sectional curvature of $(\mathcal{L}^{n\times p}_{*},\tilde{g})$ that we will need later. We use the same frame for $T_{N}\mathcal{L}^{n\times p}_{*}$ as in \cite{lin2019riemannian}: $\partial_{ii} = N_{ii}\mathcal{E}_{ii}$ and $\partial_{ij} = \mathcal{E}_{ij}, i>j$, where $\mathcal{E}_{ij}$ is an $n\times p$ matrix whose $ij$th element equals to $1$ and has zeros everywhere else. With this frame, one has,
\begin{equation*}
    \tilde{g}(\partial_{ij},\partial_{kl})=0,
\end{equation*}
which readily implies that the Christoffel symbols associated with the Riemannian product $\tilde{g}$ are all zeros \cite{absil2009optimization}:
\begin{equation} \label{flat}
    \Gamma^{ij}_{kl,uv}\equiv 0, \text{ for all } i\geq j, k\geq l, u\geq v.
\end{equation}
Be aware that we use double indices here instead of single ones as in \cite{lin2019riemannian}, but this is a trivial difference. The ``flatness" indicated by \ref{flat} will be used in subsection \ref{subsection 3.3} when we define a notion of mean.\\

\noindent
From \ref{Rie}, \ref{geo}, one can see that if $\tilde{\gamma}_{N,X}(1)=K\in\mathcal{L}^{n\times p}_{*}$, then 
\begin{equation} \label{tanvec} 
    \begin{split}
       X_{ii}/N_{ii} &= \log K_{ii}-\log N_{ii},\,\, i=1,\cdots,p\\
       X_{ij} &= K_{ij}-N_{ij},\,\, i>j, 
    \end{split}
\end{equation}
and hence the geodesic distance $\tilde{d}$ between $N, K$ is,
\begin{equation} \label{distL}
    \tilde{d}(N,K)^2 = \tilde{g}_{N}(X,X)=\sum_{i>j} (K_{ij}-N_{ij})^2+\sum_{i=1}^{p} (\log K_{ii}-\log N_{ii})^2.
\end{equation}

\begin{remark}
There were multiple purposes for the diagonal-log growth used in the coordinate chart $x$ in \cite{lin2019riemannian}: it is to prevent the swelling phenomenon in $P(n)$ (see also {\it Remark} \ref{rem:Fr} below) and to create an infinite distance between a matrix of full rank in $P(n)$ and a matrix of low rank in $S(n,p)$, for some $p<n$. There is no swelling phenomenon in $S(n,p)$ as $det(M)=0$ if $M\in S(n,p)$. We choose this diagonal-log set-up for $\mathcal{L}^{n\times p}_{*}$ and $S(n,p)^{*}$ to take advantage of the convenient formulas that follow.
\end{remark}

\subsection{The Riemannian manifold \texorpdfstring{$S(n,p)^*$}{S(n,p)^*} } \label{subsection 3.2}

\noindent
We now define a Riemannian structure on $S(n,p)^{*}$ using $\Psi$. First, we need a differential map $D_{N}\Psi: T_{N}\mathcal{L}^{n\times p}_{*}\to T_{NN^{T}}S(n,p)^{*}$. Consider the curve $\tilde{\gamma}_{N,X}(t)$ in \ref{geo}. Then 
\begin{equation} \label{geocurv} 
    S(n,p)^{*}\ni \Psi(\tilde{\gamma}_{N,X}(t))=\tilde{\gamma}_{N,X}(t)\tilde{\gamma}_{N,X}(t)^{T}
\end{equation}
for all $t\in (-\infty,\infty)$ and is a curve going through $NN^{T}\in S(n,p)^{*}$ at $t=0$; moreover,
\begin{multline} \label{res1} 
    T_{NN^{T}}S(n,p)^{*}\ni\frac{d}{dt}\Psi(\tilde{\gamma}_{N,X}(t))|_{t=0}=\tilde{\gamma}_{N,X}(0)\frac{d}{dt}\tilde{\gamma}_{N,X}(t)^{T}|_{t=0}+\frac{d}{dt}\tilde{\gamma}_{N,X}(t)|_{t=0}\tilde{\gamma}_{N,X}(0)^{T} \\=XN^{T}+NX^{T}.
\end{multline}

\noindent
We need another map, $(D_{N}\Psi)^{-1}: T_{M}S(n,p)^{*}\to T_{N}\mathcal{L}_{*}^{n\times p}$, for $NN^{T}=M$. First observe that if $\beta(t)\in S(n,p)^{*}$ is a smooth curve for $t\in (-\epsilon,\epsilon)$, then for any time $t_0$ in that interval, 
\begin{equation*} 
    S(n,p)^{*}\ni \beta(t_0)=N(t_0)N(t_0)^{T}
\end{equation*}
for a unique $N(t_0)\in\mathcal{L}^{n\times p}_{*}$ (again, by the discussion in subsection \ref{subsection 9.3}). Hence the smooth curve $\beta$ determines another smooth curve $\alpha$ such that $\alpha(t)\in\mathcal{L}^{n\times p}_{*}$ for the same time interval. Differentiating $\beta(t)=\alpha(t)\alpha(t)^{T}$ at $t=0$ gives the same formula as in \ref{res1}. In other words, any $W\in T_{M}S(n,p)^{*}$ can be written as:
\begin{equation} \label{res2} 
    W = XN^{T}+NX^{T}
\end{equation}
for some $X\in T_{N}\mathcal{L}^{n\times p}_{*}\cong\mathcal{L}^{n\times p}$. Now \ref{res2} also implies that such $X$ has to be unique. Indeed, suppose that $X$ is in the kernel of the linear map:
\begin{equation} \label{Lmap} 
    \mathcal{L}^{n\times p}\ni X\mapsto XN^{T}+NX^{T}.
\end{equation}
Then $NX^{T}=-XN^{T}=-(NX^{T})^{T}$, which implies that $NX^{T}\in\mathbb{R}^{n\times n}$ is a skew symmetric matrix. This implies that $X=0_{n\times p}$ according to the calculations in subsection \ref{subsection 9.4}. Note also that the observed surjectivity of the linear map in \ref{Lmap} alone is enough to guarantee its injectivity, as it's a linear map from a vector space of dimension $np-p(p-1)/2$ to another of the same dimension.\\

\noindent Now we augment both $N$ and $X$ in the following manners. We augment $N$ on the right with the canonical basis vectors $e_{p+1},\cdots,e_{n}$, so that the augmented matrix, called $N^{au}$, is an invertible matrix. For example, if $N=\begin{pmatrix}a_1&0\\ \star& a_2\\ \star &\star\end{pmatrix}$ for some $a_1,a_2>0,$ then $N^{au}=\begin{pmatrix}a_1&0&0\\ \star& a_2&0\\ \star&\star&1\end{pmatrix}$. We augment $X$ by adding $n-p$ zero columns to the right of it. For example, if $X=\begin{pmatrix}\star&0\\ \star&\star\\ \star &\star\end{pmatrix}$ then the augmented matrix, called $X^{au,0}$, is, $X^{au,0}=\begin{pmatrix}\star&0&0\\ \star&\star&0\\ \star&\star& 0\end{pmatrix}$. It's easily checked that, 
\begin{equation} \label{exteq} 
    Sym(n)\ni W = XN^{T}+NX^{T}=X^{au,0}(N^{au})^{T}+N^{au}(X^{au,0})^{T}.
\end{equation}
As $N^{au}$ is now a full-rank Cholesky factor of the symmetric positive definite matrix $N^{au}(N^{au})^{T}$, \ref{exteq} has a unique solution \cite{lin2019riemannian}:
\begin{equation} \label{res3} 
    X^{au,0}=N^{au}[(N^{au})^{-1}W(N^{au})^{-T}]_{1/2},
\end{equation}
where if $M\in Sym(k)$ then $[M]_{1/2}$ is the $k\times k$ lower triangular matrix obtained by keeping the lower half of $M$ and halving the diagonal elements $M_{ii}$'s. By the uniqueness discussion, we have that the $X^{au,0}$ obtained in \ref{res3} is precisely the augmented $X^{au,0}$ from $X$ in \ref{exteq}. Let $[M]^{dr}$ denote the matrix obtained from {\it dropping} the last $n-p$ columns from $M$. Then the solution $X$ for \ref{res2} is,
\begin{equation*}
    X = [N^{au}[(N^{au})^{-1}W(N^{au})^{-T}]_{1/2}]^{dr}.
\end{equation*}

\vspace{3pt}

\begin{remark}
It should be noted from \ref{res2} that not every symmetric matrix $W$ is a tangent vector in $T_{M}S(n,p)^{*}$, as $dim(Sym(n))=n(n+1)/2> np-p(p-1)/2$, only those that can be written as in \ref{res2}. That means, \ref{res3} can be used to determined whether a $W\in Sym(n)$ is in $T_{M}S(n,p)^{*}$: if $X^{au,0}$ doesn't have its last $n-p$ columns to be zeros, then $W\not\in T_{M}S(n,p)^{*}$. See also the calculations in \ref{subsection 5.3}. On the other hand, since $S(n,p)^{*}$ is an open submanifold of $S(n,p)$ (see sub-subsection \ref{subsubsection 9.3.1}), $T_{M}S(n,p)^{*}=T_{M}S(n,p)$.
\end{remark}

\subsubsection{A Riemannian choice} \label{subsubsection 3.2.1}

\noindent We now have, $(D_{N}\Psi)^{-1}: T_{M}S(n,p)^{*}\to T_{N}\mathcal{L}^{n\times p}_{*}$, for $NN^{T}=M$, with,
\begin{equation} \label{indiff} 
    (D_{N}\Psi)^{-1}: W\mapsto [N^{au}[(N^{au})^{-1}W(N^{au})^{-T}]_{1/2}]^{dr},
\end{equation}
where the notations are as in the previous subsection. We define a Riemannian metric $g$ on $S(n,p)^{*}$ via a pull-back. Let $W,V\in T_{M}S(n,p)^{*}$. Then,
\begin{equation*} 
    g_{M}(W,V) = \tilde{g}_{N}([N^{au}[(N^{au})^{-1}W(N^{au})^{-T}]_{1/2}]^{dr},[N^{au}[(N^{au})^{-1}V(N^{au})^{-T}]_{1/2}]^{dr})
\end{equation*}
or equivalently,
\begin{equation} \label{pull} 
    \tilde{g}_{N}(X,Y) = g_{\Psi(N)}(D_{N}\Psi(X),D_{N}\Psi(Y))=g_{NN^{T}}(XN^{T}+NX^{T},YN^{T}+NY^{T}).
\end{equation}
This definition \ref{pull} allows $\Psi$ to be a Riemannian isometry between $(\mathcal{L}^{n\times p}_{*},\tilde{g})$ and $(S(n,p)^{*},g)$ \cite{lee2009manifolds}. This means, for instance, the geodesics in $S(n,p)^{*}$ can be obtained from those of $\mathcal{L}^{n\times p}_{*}$ via $\Psi$; precisely, a geodesic in $S(n,p)^{*}$ has the form as in \ref{geocurv}. Hence, $S(n,p)^{*}$ is also geodesically complete, and the geodesic distance $d$ between $P,Q\in S(n,p)^{*}$, is
\begin{equation} \label{Sdist} 
    d(P,Q)=\tilde{d}(\Psi^{-1}(P),\Psi^{-1}(Q)).
\end{equation}

\subsection{A definition of mean} \label{subsection 3.3}

\noindent
Let $A_{l}=N_{l}N_{l}^{T}\in S(n,p)^{*}, l=1,\cdots,k$, we introduce the following notions of mean \begin{equation} \label{mean1} 
    \mathcal{F}(A_1,\cdots,A_{k})=\mathcal{F}(N_1,\cdots,N_{k})\mathcal{F}(N_1,\cdots,N_{k})^{T}
\end{equation}
where, $\mathcal{F}(N_1,\cdots,N_{k})\in\mathcal{L}^{n\times p}_{*}$ and
\begin{equation} \label{mean2}
    \begin{split}
        \mathcal{F}(N_1,\cdots,N_{k})_{ii} &=\prod_{l=1}^{k}\sqrt[k]{(N_{l})_{ii}},\,\,i=1,\cdots,p\\
        \mathcal{F}(N_1,\cdots,N_{k})_{ij} &= \frac{1}{k}\sum_{l=1}^{k} (N_{l})_{ij},\,\, i>j.
    \end{split}   
\end{equation}
In other words, $\mathcal{F}(N_1,\cdots,N_{k})$ registers a geometric mean along the mock diagonal and an arithmetic mean everywhere else. The reader should note that there is a slight abuse of notation in the usage of $\mathcal{F}$ in \ref{mean1}, \ref{mean2}. $\mathcal{F}$ simply means ``mean" (average), and its full meaning depends on context.\\

\noindent
Since the idea here follows {\bf 4.1} in \cite{lin2019riemannian} closely, we will give a summary, as the facts were already established in the said paper. Let $Q$ be a random point on a Riemannian manifold $\mathcal{M}$ and let $F_{Q}(x) = \mathbb{E}d^2_{\mathcal{M}}(x,Q)$, with $d_{\mathcal{M}}$ denoting the geodesic distance on $\mathcal{M}$. Then the Fr\'echet mean of $Q$, denoted $\mathcal{F}(Q)$, is defined as,
\begin{equation*}
    \mathcal{F}(Q) = \argmin_{z} F_{Q}(z)
\end{equation*}
provided that $F_{Q}(z)<\infty$ for some $z$. For an introduction to Fr\'echet and other types of mean used in statistics, see \cite{pennec2006intrinsic}. The Fr\'echet mean might not exist. However, in the case of $(\mathcal{L}^{n\times p}_{*},\tilde{g})$, if $F_{Q}$ is finite at some point, then $\mathcal{F}(Q)$ exists. This is due to its simply-connectedness (no ``holes" in the geometry of $\mathcal{L}^{n\times p}_{*}\subset\mathbb{R}^{n\times p}$) and its ``flatness" (the sectional curvature is zero due to \ref{flat}); the reasoning was well laid out in the proof of \cite[Proposition 09]{lin2019riemannian}. As $\Psi$ is a Riemannian isometry (see sub-subsection \ref{subsubsection 3.2.1}), these properties also get carried over for $(S(n,p)^{*},g)$, and hence the same is true about Fr\'echet mean on $S(n,p)^{*}$. \\

\noindent
The Fr\'echet mean on $\mathcal{L}^{n\times p}_{*}$ takes the following form:
\begin{equation} \label{form1}
    \begin{split}
        (\mathcal{F}(N))_{ii} &=\exp(\mathbb{E}\log N_{ii}),\,\, i=1,\cdots,p\\
        (\mathcal{F}(N))_{ij} &= \mathbb{E}N_{ij},\,\, i>j,
    \end{split}
\end{equation}
and that on $S(n,p)^{*}$ takes the form,
\begin{equation} \label{form2} 
    \mathcal{F}(M) = \mathcal{F}(N)\mathcal{F}(N)^{T},
\end{equation}
where $M=NN^{T}$. The proof of \ref{form1}, \ref{form2} follows exactly that of \cite[Proposition 10]{lin2019riemannian}. Now if the mean in \ref{form1} is interpreted as a sample mean - from a sample $\{N_1,\cdots, N_{k}\}$ - then one gets back \ref{mean2}, \ref{mean1} from \ref{form1}, \ref{form2}, respectively.\\

\noindent
If $k=2$, then a comparison between \ref{mean2}, \ref{mean1} and \ref{geo}, \ref{tanvec}, \ref{geocurv} shows that both \ref{mean1}, \ref{mean2} are also geodesic midpoints in $S(n,p)^{*}$ and $\mathcal{L}^{n\times p}_{*}$, respectively.

\vspace{5pt}

\begin{remark}
One can also consider a generalized version 
\begin{equation*}
    \mathcal{F}_{\vec{w}}(A_1,\cdots,A_{k}) =\mathcal{F}_{\vec{w}}(N_1,\cdots,N_{k})\mathcal{F}_{\vec{w}}(N_1,\cdots,N_{k})^{T}
\end{equation*}
of \ref{mean1}, where $\vec{w}=(w_1,\cdots,w_{k}), 0<w_{l}<1, \sum_{l=1}^{k}w_{l}=1$, and,
\begin{align*}
    \mathcal{F}(N_1,\cdots,N_{k})_{ii} & =\prod_{l=1}^{k}(N_{l})_{ii}^{w_{l}},\,\,i=1,\cdots,p\\ \mathcal{F}(N_1,\cdots,N_{k})_{ij} & = \sum_{l=1}^{k}w_{l} (N_{l})_{ij},\,\, i>j.
\end{align*}
\end{remark}

\begin{remark} \label{rem:Fr}
The Fr\'echet mean on $S(n,p)^{*}$ is convenient due to its close form \ref{mean1}. It also satisfies permutation invariance:
\begin{equation*} 
    \mathcal{F}_{\vec{w}}(A_1,\cdots,A_{k})=\mathcal{F}_{\vec{w}}(\pi(A_1,\cdots,A_{k}))
\end{equation*}
where $\pi(A_1,\cdots,A_{k})$ is any permutation of $(A_1,\cdots,A_{k})$. However, as noted about \ref{mean2}, it's neither a geometric mean nor an arithmetic mean. Hence even for the $P(n)$ case ($S(n,p)^{*}$ when $n=p$), in which it's called the log-Cholesky mean \cite{lin2019riemannian}, it doesn't satisfy a lot of physical properties in the ALM list \cite{ando2004geometric}, which is after all meant for geometric means. For example, for some $A, B\in P(n)$, $O\in O(n)$, it might happen that
\begin{equation*}
    \mathcal{F}(OAO^{T},OBO^{T})\not= O\mathcal{F}(A,B)O^{T},
\end{equation*}
i.e., the designed mean doesn't satisfy congruence invariance.
\end{remark}

\vspace{5pt}

\begin{remark}
It should be noted that a notion of mean depends on interpretation and its intended usage. In \cite{bonnabel2013rank}, the authors constructed a notion of rank-preserving geometric mean for $S(n,p)$ which satisfies three desirable properties, {\it consistency with scalars, invariance to scalings and rotations, self-duality}. Our mean is a Fr\'echet mean for matrices in $S(n,p)^{*}$, with the rank-preserving property built in, and depends straightforwardly and continuously on the entries of the involved matrices.
\end{remark}

\subsection{On the use of $S(n,p)^{*}$} \label{subsection 3.4}

\subsubsection{Covariance matrices as sole objects of investigation} \label{subsubsection 3.4.1}

\noindent
It should be noted that $S(n,p)^{*}$ as a set is not closed under similarity transformations - if $M\in S(n,p)^{*}, O\in O(n)$, it's not true that $OMO^{T}\in S(n,p)^{*}$. Closure under similarity transformations is desirable if one interprets $M\in S(n,p)^{*}$ as a linear operator. However, this lack of closure in $S(n,p)^{*}$ is not at all a detrimental effect in practice. In many applications for which this $S(n,p)^{*}$ geometry is designed, matrices arise as covariance matrices, and changing bases is not feasible. For example, the functional connectivity between ROIs (regions of interest) of the human brains is calculated as a covariance matrix $M$ using fMRI (functional MRI) data \cite{shirer2012decoding}. Here, fMRI measures brain activity by detecting changes associated with blood flow - interpreted as vectors in $\mathbb{R}^{n}$ - through $D$ discrete time points. These data form $Z\in\mathbb{R}^{n\times D}$ and subsequently the covariance matrix $M=ZZ^{T}$. It's a common assumption that the underlying covariance process is driven by a low-dimensional {\it latent} space \cite{li2019modeling}. In other words, the collected $M$ is a noisy, usually full-rank, version of the ground-truth covariance matrix $\tilde{M}$ of some dimension $p\ll n$. Now, the ROIs are predefined in a brain template according to their functions and atlas structures - in other words, ``bases" - which can't be changed \cite{shen2002hammer}. Other applications include the gene co-expression networks, which are commonly used to characterize gene-gene interactions under different disease conditions \cite{kolberg2020co} or tissue types \cite{liesecke2018ranking}. Methodically, the mRNA expressions of $p$ genes for a patient are measured as a vector $X\in\mathbb{R}^{n\times p}_{*}$, and by tradition, genes are ordered by their name alphabetically for downstream analysis. These vectors $X$'s are in turn used to estimate a covariance matrix that represents a gene co-expression network, which then can be estimated from a collection of gene expression profiles from multiple patient populations. Hence, it's fitting that we adopt the point of view that the matrices themselves are objects of investigation, and not their representations of linear operators.

\subsubsection{A brief comparison with the geometry in \cite{bonnabel2010riemannian}}

\noindent
%As alluded in the discussion above, producing meaningful averages is an important task in many statistical applications. In the case of midpoints, it often requires the existence of endpoint geodesics. 
To our knowledge, \cite{bonnabel2010riemannian} is among very few works in the literature that supply closed form, endpoint, geodesic-like curves connecting points in $S(n,p)$. Their geometry is different from ours. They inferred their metrics of $S(n,p)\cong (V(n,p)\times P(p))/O(p)$ from those of the {\it structure} space $V(n,p)\times P(p)$, while we inferred our metrics of $S(n,p)^{*}\cong \mathcal{L}^{n\times p}_{*}\times\mathcal{L}^{n\times p}_{*}$ from those of its factor space $\mathcal{L}^{n\times p}_{*}$. Both geometries are geodesically complete. On the surface, the advantage of their quotient geometry over ours is that it covers the whole space $S(n,p)$. However, their proposed metrics on $S(n,p)$ are not true geodesic metrics, even when the ones on $V(n,p)\times P(p)$ are. When it comes to midpoints, in the simple case of diagonal $A,B\in S(n,p)^{*}$, our midpoint notions coincide. For example, take
\begin{equation*} 
    A=\begin{pmatrix}a_1&0&0\\0&a_2&0\\0&0&0\end{pmatrix}, B=\begin{pmatrix}b_1&0&0\\0&b_2&0\\0&0&0\end{pmatrix}\in S(3,2)^{*} \subset S(3,2).
\end{equation*}
Then the midpoint of $A,B$ in both contexts is, $\begin{pmatrix}a_1^{1/2}b_1^{1/2}&0&0\\ 0&a_2^{1/2}b_2^{1/2}&0\\ 0&0&0\end{pmatrix}$. They behave differently otherwise. Their midpoint is a geometric midpoint by design and enjoys a longer list of desirable properties than ours (see {\it Remark} \ref{rem:Fr}). Both are rank-preserving. However, as their geometry depends on the Grassman movements, there are cases that more than one ``possible" midpoint might exist for two points $A,B\in S(n,p)$ when at least one of the principal angles, $0\leq\theta_1\leq\cdots\leq\theta_{p}\leq\pi/2$, between the two subspaces $Col(A), Col(B)\in\mathbb{R}^{n}$ is equal to $\pi/2$. See section {\bf 6} in \cite{bonnabel2010riemannian}. In contrast, our formula \ref{mean1} gives a unique Fr\'echet mean for any finite collection $A_1,\cdots, A_{k}\in S(n,p)^{*}$. It should be noted that the non-uniqueness of midpoints is not necessarily physically counterintuitive - imagine the case of two poles on the sphere. Our geometry evades such non-uniqueness by being restricted. Take $A=\begin{pmatrix} 1&1\\ 1&1\end{pmatrix}, B=\begin{pmatrix} 1&-1\\ -1&1\end{pmatrix}\in S(2,1)^{*}\subset S(2,1)$. Clearly, the only principle angle is $\theta_1=\pi/2$. Their geometry doesn't give a unique midpoint between $A,B$, but ours does: $\begin{pmatrix}1&0\\ 0&0\end{pmatrix}$.

\subsubsection{$S(n,p)^{*}$ and $S(n,p)$}

\noindent
We work with $S(n,p)^{*}$ because we can exploit its convenient {\it reduced} Cholesky factorization and because $S(n,p)^{*}$ is dense in $S(n,p)$ (see sub-subsection \ref{subsubsection 9.3.3}). Now suppose that $A\in S(n,p)^{*}$ and $B\in S(n,p)-S(n,p)^{*}$. For concreteness, suppose $n=3, p=2$. Let $A = NN^{T}, B=KK^{T}$, where
\begin{equation*}
    N=\begin{pmatrix}n_{11}&0\\ n_{21}&n_{22}\\n_{31}&n_{32}\end{pmatrix},\text{ and } K=\begin{pmatrix}0&0\\ k_{21}&0\\k_{31}&k_{32}\end{pmatrix},
\end{equation*}
where $n_{11}, n_{22}>0$. Then as remarked in sub-subsection \ref{subsubsection 9.3.3}, we can approximate $B$ with $K_{\epsilon}K_{\epsilon}^{T}$, where
\begin{equation*}
    K_{\epsilon} = \begin{pmatrix}\epsilon&0\\ k_{21}&\epsilon\\ k_{31}&k_{32}\end{pmatrix}.
\end{equation*}
A straightforward calculation gives
\begin{equation} \label{ineq0} 
    \|KK^{T}-K_{\epsilon}K_{\epsilon}^{T}\|_{F}=\|B-K_{\epsilon}K_{\epsilon}^{T}\|_{F}\geq \epsilon\|K\|_{F}.
\end{equation}

\noindent
Now the midpoint $\mathcal{F}(A,K_{\epsilon}K_{\epsilon}^{T})=JJ^{T}$ with,
\begin{equation*}
    J=\begin{pmatrix}\sqrt{n_{11}}\sqrt{\epsilon}&0\\ (n_{21}+k_{21})/2 &\sqrt{n_{22}}\sqrt{\epsilon}\\ (n_{31}+k_{31})/2&(n_{32}+k_{32})/2\end{pmatrix}.
\end{equation*}
The geodesic distance between $JJ^{T}$ and $K_{\epsilon}K_{\epsilon}^{T}$ according to \ref{distL}, \ref{Sdist} is, 
\begin{equation} \label{ineq1} 
    d(JJ^{T},K_{\epsilon}K_{\epsilon}^{T})^2=\tilde{d}(J,K_{\epsilon})^2=\sum_{i>j}(n_{ij}/2-k_{ij}/2)^2+\sum_{i=1}^2((\log n_{ii})/2-(\log\epsilon)/2)^2.
\end{equation}
Identify $E\in\mathbb{R}^{3\times 2}$ with its augmented version $E^{au,0}\in\mathbb{R}^{3\times 3}$ (recall the notation from subsection \ref{subsection 3.2}). We then have the following facts,
\begin{equation*}
    \|E\|_{F}=\|E^{au,0}\|_{F}\,\,\,\text{ and }\,\,\, \|E^{au,0}\|_{F}\asymp_{p}\|E^{au,0}\|_{op},
\end{equation*}
where $\|\cdot\|_{op}$ is the induced operator norm of $E^{au,0}: (\mathbb{R}^3,\|\cdot\|_2) \to (\mathbb{R}^3,\|\cdot\|_2)$. Now if $E,H\in\mathbb{R}^{3\times 3}$ then \cite{davidson1996c}
\begin{equation*}
    \|EH\|_{op}\leq\|E\|_{op}\|H\|_{op}.
\end{equation*}
It follows that the Frobenius distance between $JJ^{T}$ and $B=KK^{T}$ is then dominated by,
\begin{multline} \label{ineq2} 
    \|JJ^{T}-KK^{T}\|_{F}\leq \|J(J-K)^{T}\|_{F}+\|(J-K)K^{T}\|_{F}\lesssim_{p}(\|J^{au,0}\|_{op}\|+\|K^{au,0}\|_{op})\|(J-K)^{au,0}\|_{op}\\ \lesssim_{p} c\max\{\|K\|_{\infty},\|J\|_{\infty}\}\|K-J\|_{F}\leq c\max\{\|K\|_{\infty},\|J\|_{\infty}\}\bigg(\sum_{i>j}(n_{ij}/2-k_{ij}/2)^2+\sum_{i=1}^2\epsilon n_{ii}\bigg)^{1/2},
\end{multline}
where $\|M\|_{\infty}=\max_{i,j}|M_{ij}|$. Observe that if one takes $\epsilon$ sufficiently small, so that,
\begin{equation} \label{cond} 
    \log\sqrt{n_{ii}}\geq\log\sqrt{\epsilon}+\sqrt{\epsilon}\sqrt{n_{ii}}
\end{equation}
then \ref{ineq0}, \ref{ineq1}, \ref{ineq2} imply that,
\begin{multline} \label{conc} 
    \|KK^{T}-K_{\epsilon}K_{\epsilon}^{T}\|_{F}+d(JJ^{T},K_{\epsilon}K_{\epsilon}^{T})\geq\epsilon\|K\|_{F}+\bigg(\sum_{i>j}(n_{ij}/2-k_{ij}/2)^2+\sum_{i=1}^2((\log n_{ii})/2-(\log\epsilon)/2)^2\bigg)^{1/2}\\  \gtrsim\max\{\|K\|_{\infty},\|J\|_{\infty}\}\bigg(\sum_{i>j}(n_{ij}/2-k_{ij}/2)^2+\sum_{i=1}^2\epsilon n_{ii}\bigg)^{1/2}\geq\|JJ^{T}-KK^{T}\|_{F}.
\end{multline}
In short, \ref{conc}, which is of triangle-inequality type, says that the usual Frobenius distance between $KK^{T}\in S(n,p)-S(n,p)^{*}$ to the midpoint $JJ^{T}$ of $NN^{T}$ and the approximation $K_{\epsilon}K_{\epsilon}^{T}$ is controllable by a constant multiple of the sum $\|KK^{T}-K_{\epsilon}K_{\epsilon}^{T}\|_{F}+d(JJ^{T},K_{\epsilon}K_{\epsilon}^{T})$ - i.e., the midpoint found through approximation cannot be far away from $KK^{T}$. The implicit constant in \ref{conc} is dependent on $\max\{\|K\|_{\infty},\|J\|_{\infty}\}$ and on $p$ for general $n,p$. Moreover, for general $n,p$, the condition \ref{cond} is still sufficient and \ref{conc} still holds, as $K$ is the worst case scenario: all $k_{ii}=0$.\\

\noindent
Suppose one encounters a rarer case where both $A=NN^{T}, B=KK^{T}\in S(n,p)-S(n,p)^{*}$, say, $A,B\in S(3,2)-S(3,2)^{*}$ and, 
\begin{equation*}
    N=\begin{pmatrix}0&0\\ n_{21}& 0\\n_{31}&n_{32}\end{pmatrix},\text{ and } K=\begin{pmatrix}0&0\\ k_{21}&0\\k_{31}&k_{32}\end{pmatrix},
\end{equation*}
then one can approximate $A, B$ with $N_{\epsilon_1}N_{\epsilon_1}^{T},K_{\epsilon_2}K_{\epsilon_2}^{T}$, respectively, where,
\begin{equation*}
    N_{\epsilon_1} = \begin{pmatrix}\epsilon_1&0\\ n_{21}&\epsilon_1\\ n_{31}&n_{32}\end{pmatrix},\text{ and } K_{\epsilon_2} = \begin{pmatrix}\epsilon_2&0\\ k_{21}&\epsilon_2\\ k_{31}&k_{32}\end{pmatrix}.
\end{equation*}
An alternative sufficient version of \ref{cond} is now,
\begin{equation*}
    |\log\sqrt{\epsilon_1}-\log\sqrt{\epsilon_2}|\geq\sqrt{\epsilon_1}\sqrt{\epsilon_2},
\end{equation*}
and the following version of \ref{conc} holds:
\begin{align*}
    \|NN^{T}-N_{\epsilon_1}N_{\epsilon_1}^{T}\|_{F}+d(JJ^{T},N_{\epsilon_1}N_{\epsilon_1}^{T}) &\gtrsim_{J,N}\|JJ^{T}-NN^{T}\|_{F},\\ 
    \|KK^{T}-K_{\epsilon_2}K_{\epsilon_2}^{T}\|_{F}+d(JJ^{T},K_{\epsilon_2}K_{\epsilon_2}^{T}) &\gtrsim_{J,K}\|JJ^{T}-KK^{T}\|_{F},
\end{align*}
where $JJ^{T}=\mathcal{F}(N_{\epsilon_1}N_{\epsilon_1}^{T},K_{\epsilon_2}K_{\epsilon_2}^{T})$.\\

\noindent
This argument, together with sub-subsections \ref{subsubsection 3.4.1}, \ref{subsubsection 9.3.3}, gives us another justification to consider the convenient geometry of $S(n,p)^{*}$.

\subsection{A geometric example}

\noindent
Let $n=2, p=1$. Then $S(2,1)^{*}$ is a manifold of dimension $2$. Moreover, since $S(n,p)$ is a conic manifold with one connected component \cite{helmke1995critical}, $S(2,1)$ can be identified with a two-dimensional cone. If $M\in S(2,1)$ then its reduced Cholesky factor $N$ is either in $\mathcal{L}^{2\times 1}_{*}$, which means $N$ is a vector in the half plane $\bigg\{\begin{pmatrix}x\\y\end{pmatrix}: x>0\bigg\}$, or $N=\begin{pmatrix}0\\y\end{pmatrix}$ with $y>0$. An illustration of $S(2,1), S(2,1)^{*}, \mathcal{L}^{2\times 1}_{*}$ is given below in Figure \ref{fig:cone}A and \ref{fig:cone}B.\\

\noindent
Note that everything that is outside of the ``slit" on the cone $S(2,1)$ belongs to $S(2,1)^{*}$. That slit is precisely the set $C_{\{2\}}$ (see the notation in section \ref{section 2}). In this case $C_{\{2\}}$ has a ray-like structure, as $\alpha\begin{pmatrix}0\\y\end{pmatrix}\begin{pmatrix}0\\y\end{pmatrix}^{T}=\alpha\begin{pmatrix}0&0\\ 0&y^2\end{pmatrix}\in C_{\{2\}}$ for $\alpha>0$ (where $y>0$). That slit is also where the Cholesky decomposition is discontinuous on $S(2,1)$. For example, $M_1=\begin{pmatrix}0&0\\ 0&1\end{pmatrix}\in S(2,1)-S(2,1)^{*}$. Then its reduced Cholesky factor is $N_1=\begin{pmatrix}0\\ 1\end{pmatrix}$. Both of the matrices,
\begin{equation*}
    M_2=\begin{pmatrix}\epsilon^2&\epsilon\\ \epsilon&1\end{pmatrix},\,\,\,\,\,\,\,\, M_3=\begin{pmatrix}\epsilon^2& -\epsilon\\ -\epsilon&1\end{pmatrix},
\end{equation*}
are close to $M_1$ in the Euclidean distance and hence are close to each other, as $\epsilon\to 0^{+}$. The reduced Cholesky factors for $M_2, M_3$ are, respectively,
\begin{equation*} 
    N_2=\begin{pmatrix}\epsilon\\1\end{pmatrix},\,\,\,\,\,\,\,\, N_3=\begin{pmatrix}\epsilon\\-1\end{pmatrix}.
\end{equation*}
These two vectors are clearly distanced away from each other: $\|N_2-N_3\|_2=\sqrt{2}$. This discontinuity is also illustrated in Figure \ref{fig:cone}, (A) and (B).

\begin{figure}
    \centering
    \includegraphics[scale=0.7]{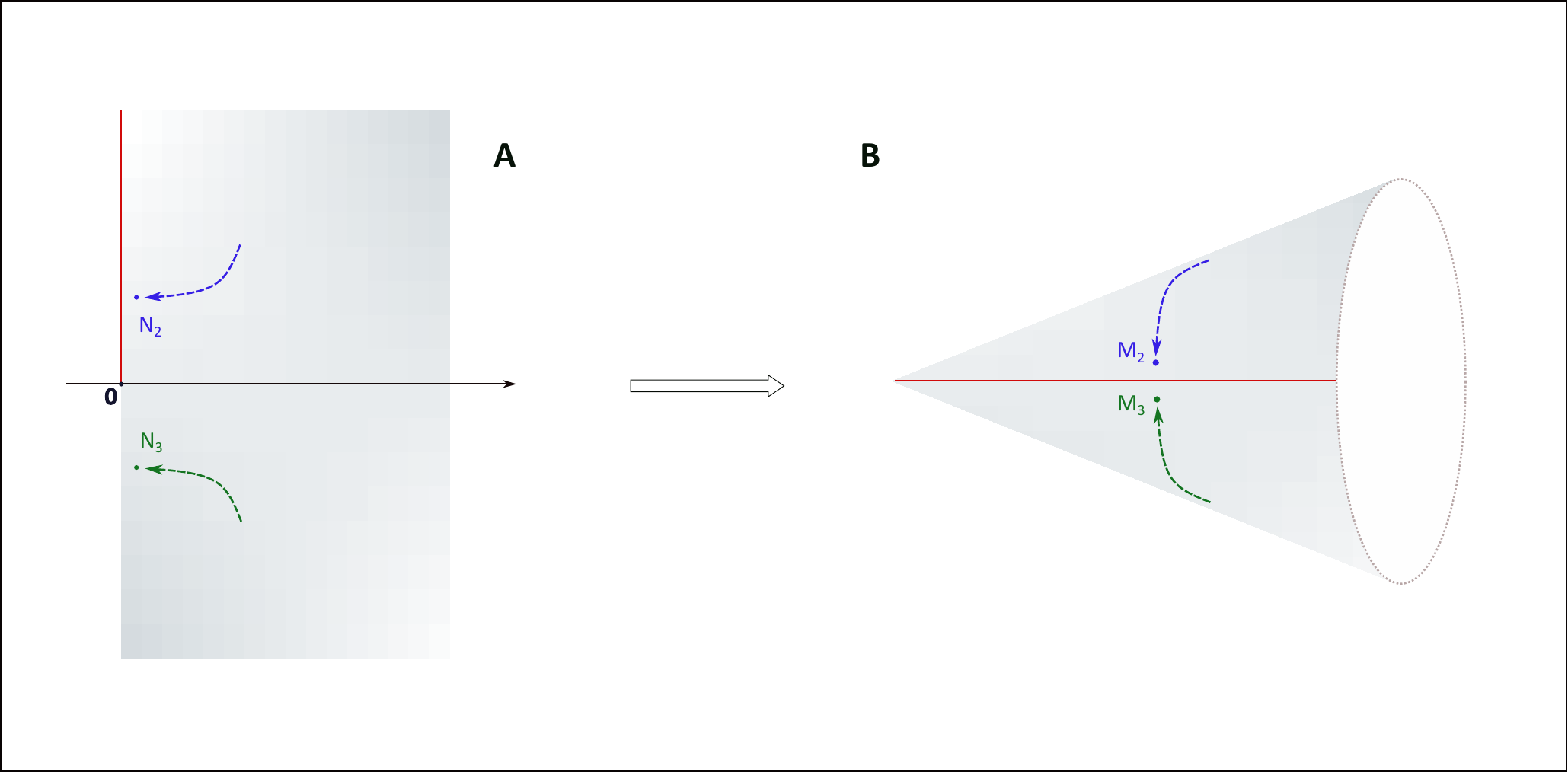}
    \caption{(A) The $xy$-plane with the shaded region representing $\mathcal{L}^{2\times 1}_{*}$; the positive $y$-axis represents the reduced Cholesky factors of matrices in $C_{\{2\}}$. \\
    (B) The cone $S(2,1)$ with the slit $C_{\{2\}}$. Two sequences converging to the same point on the slit can be sent to two divergent sequences of reduced Cholesky factors.}
    \label{fig:cone}
\end{figure}

\begin{figure}
    \centering
    \includegraphics[scale=1.15]{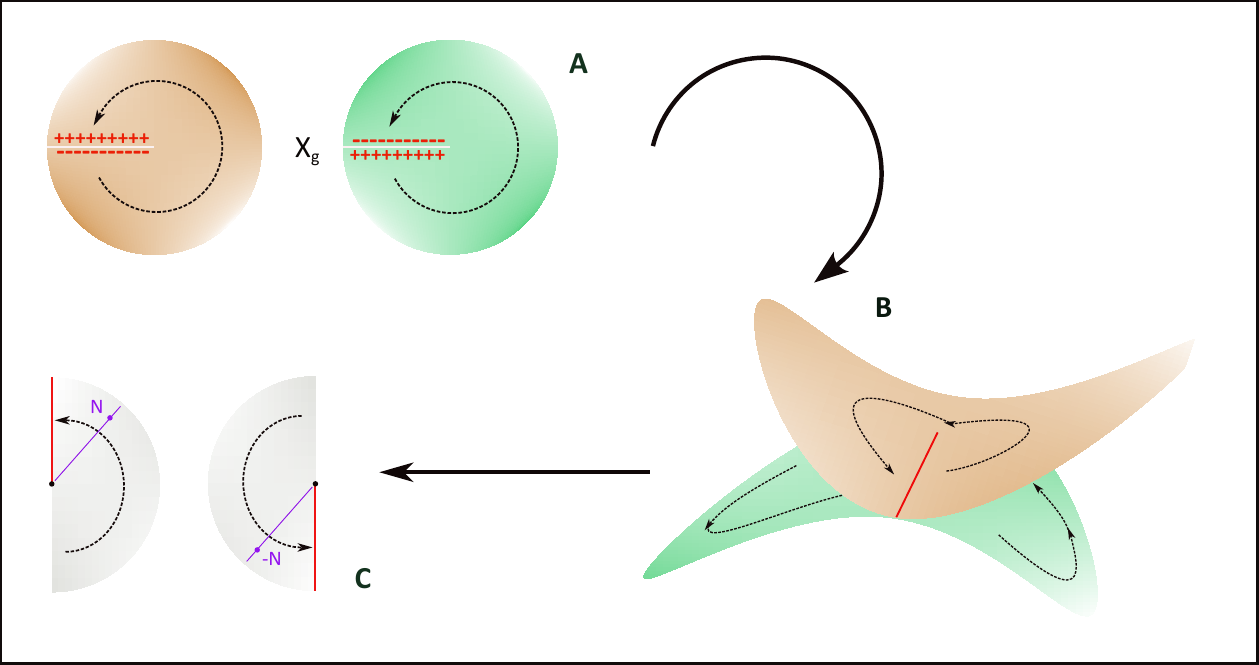}
    \caption{(A) The glued structure $S(2,1)\times_{g} S(2,1)$: the two different colored surfaces are the two copies of $S(2,1)$, each identified with $\mathbb{R}^2-\{\vec{0}\}$. A point $M\in S(2,1)$ is shown traveling from one copy to another via the glued slit.\\
    (B) The cone $S(2,1)$ is identified with $\mathbb{R}^2-\{\vec{0}\}$; the dark red line represents the slit $C_{\{2\}}$. The dashed line shows a point traveling throughout the glued structure $S(2,1)\times_{g} S(2,1)$.\\
    (C) Two halves of $\mathbb{R}^2-\{0\}$; as $M$ travels throughout the top copy of $S(2,1)$ in (A), $\psi$ outputs the reduced Cholesky factor of $M$, and when $M$ reaches the bottom copy, $\psi$ outputs the negative of its reduced Cholesky factor.}
    \label{fig:Cholesky_map_generalized}
\end{figure}

\vspace{5pt}

\begin{remark}
The picture given in Figure \ref{fig:cone} is representative for general $n,p$, since $S(n,p)$ still has a conic structure, and if $M\in C_{\{i_1,\cdots,i_{p}\}}\subset S(n,p)$ then $\alpha M\in C_{\{i_1,\cdots,i_{p}\}}$, for $\alpha>0$. One still also has that $S(n,p)^{*}$ is path-connected and dense in $S(n,p)=S(n,p)^{*}\cup\bigg(\bigcup_{\{i_1,\cdots,i_{p}\}}C_{\{i_1,\cdots,i_{p}\}}\bigg)$.
\end{remark}

\vspace{5pt}

\begin{remark}
We note in passing that one can use an elementary complex analysis trick \cite{gamelin2003complex} to define a continuous Cholesky-like map on a special ``glued" structure $S(2,1)\times_{g} S(2,1)$ as follows. Since $S(2,1)$ is a smooth two-dimensional cone, one can ``flatten" it out and identify it with $\mathbb{R}^2-\{\vec{0}\}$ - a rigorous proof of this identification can be done via applying the Uniformization Theorem for Riemannian Surfaces \cite{gamelin2003complex}. The structure $S(2,1)\times_{g} S(2,1)$ is formed by gluing the top of the slit $C_{\{2\}}$ of the first copy of $S(2,1)$ to the bottom of the slit $C_{\{2\}}$ of the second copy, and similarly, gluing the bottom of the slit $C_{\{2\}}$ of the first copy to the top of the slit $C_{\{2\}}$ of the second copy. This ``gluing" (identification) is done so that only points near the slits that are glued to each other are considered ``close". That means, for instance, points near the top slit and points near the bottom slit of the first copy are not close to each other. See Figure \ref{fig:Cholesky_map_generalized}A.
\end{remark}

\noindent
Let $\psi: S(2,1)\times_{g} S(2,1)\to\mathbb{R}^2-\{\vec{0}\}$ be as follows. On the first copy of $S(2,1)$, $\psi$ behaves like the reduced Cholesky map. On the second copy, $\psi$ is the negation of the reduced Cholesky map. Figure \ref{fig:Cholesky_map_generalized} shows $\psi$ in action as $M$ traverses throughout $S(2,1)\times_{g} S(2,1)$. It's easy to see that $\psi$ is a continuous map on $S(2,1)\times_{g} S(2,1)$. This is precisely the technique used to construct the Riemannian domain for the complex-valued square root map \cite{gamelin2003complex}.

\section{Exponentials and logarithms} \label{section 4}

\noindent As is the custom for a Riemmanian geometry, we now define the exponential and logarithm maps on $S(n,p)^{*}$. Similar to the previous section, we do so by first defining the same maps on $\mathcal{L}^{n\times p}_{*}$ and then designing the corresponding ones on $S(n,p)^{*}$ using the map $\Psi$ in \ref{psi1}.\\

\noindent
Recall the structure of geodesic curves on $\mathcal{L}^{n\times p}_{*}$ in \ref{geo} and the fact that $\mathcal{L}^{n\times p}_{*}$ is geodesically complete. Then by definition \cite{lee2009manifolds}, $\widetilde{Exp}_{N}(X)=\tilde{\gamma}_{N,X}(1)$. The logarithm map at $N$, denoted $\widetilde{Log}_{N}$, is easily defined, as it is the inverse of $\widetilde{Exp}_{N}$. That means, $\widetilde{Log}_{N}(K)=X$ if $\widetilde{Exp}_{N}(X)=K$. In short, we have, 
\begin{equation*} 
    \widetilde{Exp}_{N}: T_{N}\mathcal{L}^{n\times p}_{*}\cong\mathcal{L}^{n\times p}\ni X\mapsto \widetilde{Exp}_{N}(X)\in\mathcal{L}^{n\times p}_{*},
\end{equation*}
where, as from \ref{geo},
\begin{align*} 
    (\widetilde{Exp}_{N}(X))_{ii} &= N_{ii}\exp(X_{ii}N_{ii}^{-1}),\,\, i=1,\cdots,p\\
    (\widetilde{Exp}_{N}(X))_{ij} &= N_{ij}+X_{ij},\,\, i>j,
\end{align*}
and, 
\begin{equation*}
    \widetilde{Log}_{N}: \mathcal{L}^{n\times p}_{*}\ni K\mapsto\widetilde{Log}_{N}(K)\in T_{N}\mathcal{L}^{n\times p}_{*}\cong\mathcal{L}^{n\times p},
\end{equation*}
with, 
\begin{align*} 
    (\widetilde{Log}_{N}(K))_{ii} &= (\log K_{ii}-\log N_{ii})N_{ii},\,\, i=1,\cdots,p\\
    (\widetilde{Log}_{N}(K))_{ij} &= K_{ij}-N_{ij},\,\, i>j,
\end{align*}
as noted in \ref{tanvec}.\\

\noindent
As discussed in sub-subsection \ref{subsubsection 3.2.1}, the geodesic curves on $S(n,p)^{*}$ take the form of $\tilde{\gamma}_{N,X}(t)\tilde{\gamma}_{N,X}(t)^{T}$ (\ref{geocurv}), for some $N\in\mathcal{L}^{n\times p}_{*}$ and $X\in T_{N}\mathcal{L}^{n\times p}_{*}$. Using the differential maps $D_{N}\Psi$ and their corresponding inverses $(D_{N}\Psi)^{-1}$ we conclude that the geodesic curves $\gamma_{M,W}(t)$ on $S(n,p)^{*}$, with $M\in S(n,p)^{*}$ and $W\in T_{M}S(n,p)^{*}$, are as follows,
\begin{equation*} 
    \gamma_{M,W}(t)=\big(\tilde{\gamma}_{\Psi^{-1}(M),(D_{N}\Psi)^{-1}(W)}(t)\big)\big(\tilde{\gamma}_{\Psi^{-1}(M),(D_{N}\Psi)^{-1}(W)}(t)\big)^{T},
\end{equation*}
where $NN^{T}=M$ (or $\Psi^{-1}(M)=N$) and $(D_{N}\Psi)^{-1}$ is defined in \ref{indiff}. Just as $\mathcal{L}^{n\times p}_{*}$, $S(n,p)^{*}$ is geodesically complete (see sub-subsection \ref{subsubsection 3.2.1}). Hence the exponential map $Exp_{M}$ on $S(n,p)^{*}$ at $M$ is, 
\begin{equation*}
    Exp_{M}: T_{M}S(n,p)^{*}\ni W\mapsto Exp_{M}(W)\in S(n,p)^{*},
\end{equation*}
where, 
\begin{align} \label{exp} 
    Exp_{M}(W) = \gamma_{M,W}(1)  &=\big(\tilde{\gamma}_{\Psi^{-1}(M),(D_{N}\Psi)^{-1}(W)}(1)\big)\big(\tilde{\gamma}_{\Psi^{-1}(M),(D_{N}\Psi)^{-1}(W)}(1)\big)^{T}\\ 
    \nonumber &= \big(\widetilde{Exp}_{N}([N^{au}[(N^{au})^{-1}W(N^{au})^{-T}]_{1/2}]^{dr})\big)\big(\widetilde{Exp}_{N}([N^{au}[(N^{au})^{-1}W(N^{au})^{-T}]_{1/2}]^{dr})\big)^{T}.
\end{align}
Finally, the logarithm map $Log_{M}$ on $S(n,p)^{*}$ is as follows,
\begin{equation} \label{log} 
    Log_{M}: S(n,p)^{*}\ni Q\mapsto Log_{M}(Q)\in T_{M}S(n,p)^{*},
\end{equation}
and,
\begin{equation*} 
    Log_{M}(Q) = (D_{\Psi^{-1}(M)}\Psi)(\widetilde{Log}_{\Psi^{-1}(M)}\Psi^{-1}(Q)) = (D_{N}\Psi)(\widetilde{Log}_{N}\Psi^{-1}(Q)).
\end{equation*}
To see the last one, note that if $Log_{M}(Q)=W$ then $Exp_{M}(W)=Q$. Let 
\begin{equation*} 
    X = \widetilde{Log}_{\Psi^{-1}(M)}\Psi^{-1}(Q) = \widetilde{Log}_{N}\Psi^{-1}(Q) \text{ and } W =(D_{N}\Psi)(\widetilde{Log}_{N}\Psi^{-1}(Q)).
\end{equation*} 
Then from \ref{indiff}, $X = [N^{au}[(N^{au})^{-1}W(N^{au})^{-T}]_{1/2}]^{dr}$ and
\begin{align*}
    Exp_{M}(W) &= (\widetilde{Exp}_{N}X)(\widetilde{Exp}_{N}X)^{T} \\ 
    &=(\widetilde{Exp}_{N}(\widetilde{Log}_{N}\Psi^{-1}(Q)))(\widetilde{Exp}_{N}(\widetilde{Log}_{N}\Psi^{-1}(Q)))^{T}\\&=\Psi^{-1}(Q)\Psi^{-1}(Q)^{T}=Q,
\end{align*}
as expected.

\section{Parallel transport} \label{section 5}

\noindent In this section, we set out to find a formulation for parallel transport of tangent vectors along geodesics on $S(n,p)^{*}$; such a tool is essential in some applications of statistical analysis or machine learning. The task is made simple once we can prove that $S(n,p)^{*}$ possesses an abelian Lie group structure in which $g$ is a bi-variant metric, i.e., $g$ is invariant with respect to the left and right translations induced by the group action, and invoke the following lemma from \cite{lin2019riemannian}. 

\begin{lemma} \label{lem:lin}
Let $(\mathcal{G},\cdot)$ be an abelian Lie group with a bi-variant metric. The parallel transporter $\tau_{p,q}$ that transports tangent vectors at $p$ to tangent vectors at $q$ along the geodesic connecting $p, q$ is given by 
\begin{equation*}
    \tau_{p,q}(u)=(D_{p}l_{q\cdot p^{-1}})u
\end{equation*}
for $u\in T_{p}\mathcal{G}$, where $l_{q\cdot p^{-1}}$ stands for the left translation by $q\cdot p^{-1}\in\mathcal{G}$.
\end{lemma}

\noindent In subsection \ref{subsection 5.1} below, we prepare a Lie group structure for $S(n,p)^{*}$, in order to apply {\bf Lemma} \ref{lem:lin} in subsection \ref{subsection 5.2}. The final formulation for parallel transport is given in \ref{par}. In the last subsection \ref{subsection 5.3}, we give a sample calculation that illustrates the use of \ref{par}. Since the example there is presented in a much simplified manner (see \ref{parex}), the impatient reader might want to skip ahead to that subsection. 

\subsection{A Lie group structure} \label{subsection 5.1}

\noindent
Define a binary operation $\oplus$ on $\mathcal{L}^{n\times p}_{*}$,
\begin{align*} 
    (N_1 \oplus N_2)_{ii} &= (N_1)_{ii}\cdot (N_2)_{ii},\,\,i=1,\cdots,p,\\
    (N_1 \oplus N_2)_{ij} &= (N_1)_{ij}+ (N_2)_{ij},\,\, i>j.
\end{align*}
This operation makes $\mathcal{L}^{n\times p}_{*}$ into an abelian group. It's clear that the group identity is $I_{n\times p}$, where $(I_{n\times p})_{ii} = 1,\,\, i=1,\cdots,p,$ and $(I_{n\times p})_{ij} = 0,\,\, i\not=j$, while the group inverse $N^{-1}$ of $N\in\mathcal{L}^{n\times p}_{*}$ is,
\begin{align*} 
    (N^{-1})_{ii} &= N_{ii}^{-1},\,\, i=1,\cdots,p,\\
    (N^{-1})_{ij} &= -N_{ij},\,\, i>j.
\end{align*}
The left translation by $K\in\mathcal{L}^{n\times p}_{*}$, is $K\oplus\cdot: N\mapsto K\oplus N$. Let $\tilde{\gamma}_{N,X}(t)$ be the geodesic curve starting at $N\in\mathcal{L}^{n\times p}_{*}$ and following $X\in T_{N}\mathcal{L}^{n\times p}_{*}\cong\mathcal{L}^{n\times p}$ as in \ref{geo}. Then, $(K\oplus\tilde{\gamma}_{N,X}(t))_{ii} = K_{ii}N_{ii}\exp(tX_{ii}N_{ii}^{-1})$ and $(K\oplus\tilde{\gamma}_{N,X}(t))_{ij} = K_{ij}+N_{ij}+tX_{ij}, i>j$. Taking $\frac{d}{dt}(K\oplus\tilde{\gamma}_{N,X}(t))|_{t=0}$, we have,
\begin{equation} \label{deriv} 
    D_{N}(K\oplus\cdot): T_{N}\mathcal{L}^{n\times p}_{*}\cong\mathcal{L}^{n\times p}\ni X\mapsto D_{N}(K\oplus\cdot)(X)=\frac{d}{dt}(K\oplus\tilde{\gamma}_{N,X}(t))|_{t=0}\in T_{K\oplus N}\mathcal{L}^{n\times p}_{*}\cong\mathcal{L}^{n\times p},
\end{equation}
where the arrival vector is such that 
\begin{equation} \label{dis}
    \begin{split}
        (D_{N}(K\oplus\cdot)(X))_{ii} &= K_{ii}X_{ii},\,\, i=1,\cdots,p\\ 
        (D_{N}(K\oplus\cdot)(X))_{ij} &= X_{ij},\,\, i>j.
    \end{split}  
\end{equation}
Note that derivative map $D_{N}(K\oplus\cdot)$ doesn't depend on the location $N$ and that 
\begin{equation*} 
    D_{N}(K^{-1}\oplus\cdot)\circ D_{N}(K\oplus\cdot) = Id|_{T_{N}\mathcal{L}^{n\times p}_{*}}.
\end{equation*}  
As the operation $\oplus$ clearly commutes, the right translation by $K$, $N\mapsto N\oplus K$, is the same as the left translation by $K$, and $D_{N}(K\oplus\cdot)=D_{N}(\cdot\oplus K)$. All of this makes $\mathcal{L}^{n\times p}_{*}$ a Lie group \cite{helgason1979differential}. Moreover, it follows from \ref{Rie} that
\begin{multline} \label{bi} 
    \tilde{g}_{N\oplus K}(D_{N}(\cdot\oplus K)(X),D_{N}(\cdot\oplus K)(Y))=\tilde{g}_{K\oplus N}(D_{N}(K\oplus\cdot)(X),D_{N}(K\oplus\cdot)(Y)) \\=\sum_{i>j} X_{ij}Y_{ij}+\sum_{i=1}^{p} N_{ii}^{-2}X_{ii}Y_{ii} = \tilde{g}_{N}(X,Y),
\end{multline}
for $X,Y\in T_{N}\mathcal{L}^{n\times p}_{*}$ and $D_{N}(\cdot\oplus K)(X),D_{N}(\cdot\oplus K)(Y), D_{N}(K\oplus\cdot)(X),D_{N}(K\oplus\cdot)(Y)\in T_{K\oplus N}\mathcal{L}^{n\times p}_{*}$ (see \ref{deriv}). Now \ref{bi} says that $\tilde{g}$ is a bi-variant metric on $(\mathcal{L}^{n\times p}_{*},\oplus)$. We want to have the same thing for $S(n,p)^{*}$.\\

\noindent We define a binary operation $\otimes$ on $S(n,p)^{*}$ as follows:
\begin{equation} \label{otimes} 
    P_1 \otimes P_2 = \Psi(\Psi^{-1}(P_1)\oplus\Psi^{-1}(P_2))= (\Psi^{-1}(P_1)\oplus\Psi^{-1}(P_2))(\Psi^{-1}(P_1)\oplus\Psi^{-1}(P_2))^{T}.
\end{equation}
It's clear that like $\oplus$, $\otimes$ is also an abelian operation. Moreover, the group inverse of $P\in S(n,p)^{*}$ is, $(\Psi^{-1}(P))^{-1}((\Psi^{-1}(P))^{-1})^{T}=\Psi((\Psi^{-1}(P))^{-1})$. The identity for $\otimes$ is $\Psi(I_{n\times p})$. While associativity is clear for $\oplus$, the same for $\otimes$ needs some checking:
\begin{align*} 
    (P_1 \otimes P_2)\otimes P_3 &=\Psi(\Psi^{-1}(P_1)\oplus\Psi^{-1}(P_2))\otimes P_3\\
    &=\Psi\big((\Psi^{-1}(P_1)\oplus\Psi^{-1}(P_2))\oplus\Psi^{-1}(P_3)\big)\\
    &=\Psi\big(\Psi^{-1}(P_1)\oplus(\Psi^{-1}(P_2)\oplus\Psi^{-1}(P_3))\big)\\
    &=P_1\otimes(P_2\otimes P_3).
\end{align*}
Hence, 
\begin{equation*}
    \Psi: (\mathcal{L}^{n\times p}_{*},\oplus)\to (S(n,p)^{*},\otimes)
\end{equation*}
is a group isomorphism \cite{dummit2004abstract}. Recall from the definition \ref{pull} that $\Psi$ is also a Riemannian isometry between $(\mathcal{L}^{n\times p}_{*},\tilde{g})$ and $(S(n,p)^{*},g)$. Now the left and right translations on $S(n,p)^{*}$, with respect to $\otimes$, as well as their derivatives, can be defined similarly to the case in $(\mathcal{L}^{n\times p}_{*},\oplus)$. For instance, following \ref{deriv}, the derivative map of the left translation by $Q\in S(n,p)^{*}$ at $P\in S(n,p)^{*}$ maps:
\begin{multline} \label{deriv2} 
    T_{P}S(n,p)^{*}\ni W\mapsto \frac{d}{dt}\Psi(N_{Q}\oplus\tilde{\gamma}_{N_{P},X}(t))|_{t=0} \\=(D_{N_{P}}(N_{Q}\oplus\cdot)(X))(N_{Q}\oplus N_{P})^{T}+(N_{Q}\oplus N_{P})(D_{N_{P}}(N_{Q}\oplus\cdot)(X))^{T}\in T_{Q\otimes P}S(n,p)^{*},
\end{multline}
where $N_{P}=\Psi^{-1}(P), N_{Q}=\Psi^{-1}(Q)$ and $X=(D_{N_{P}}\Psi)^{-1}(W)=[N^{au}_{P}[(N^{au})^{-1}_{P}W(N^{au})^{-T}_{P}]_{1/2}]^{dr}$ defined in \ref{indiff}. \\

\noindent More can be said from \ref{deriv2}. Let $N=N_{Q}\oplus N_{P}$. Note then from \ref{otimes}, $\Psi^{-1}(Q\otimes P)=\Psi^{-1}(Q)\oplus\Psi^{-1}(P)=N$. Recall that the differential map $D_{N}\Psi$ in \ref{res1} maps:
\begin{equation*} 
    T_{N}\mathcal{L}^{n\times p}_{*}\ni Z\mapsto ZN^{T}+NZ^{T}\in T_{\Psi(N)}S(n,p)^{*}=T_{Q\otimes P}S(n,p)^{*}.
\end{equation*}
Hence by matching \ref{deriv2} with \ref{res1}, one realizes that \ref{deriv2} also says,
\begin{multline} \label{obs} 
    D_{N}\Psi: D_{N_{P}}(N_{Q}\oplus\cdot)(X)\mapsto D_{N}\Psi(D_{N_{P}}(N_{Q}\oplus\cdot)(X))=(D_{N_{P}}(N_{Q}\oplus\cdot)(X))N^{T}+N(D_{N_{P}}(N_{Q}\oplus\cdot)(X))^{T}\\ =D_{P}(Q\otimes\cdot)(W).
\end{multline}

\noindent Now it's easy to show that, through $\Psi$, $g$ is a bi-variant metric on $S(n,p)^{*}$, as $\tilde{g}$ is on $\mathcal{L}^{n\times p}_{*}$. We need to prepare a few notations first. Let $X=(D_{N_{P}}\Psi)^{-1}(W)$ and $Y=(D_{N_{P}}\Psi)^{-1}(V)$. From \ref{obs}, we have
\begin{align} \label{X} 
    D_{N}\Psi(D_{N_{P}}(N_{Q}\oplus\cdot)(X)) &= D_{P}(Q\otimes\cdot)(W)\\
    \label{Y} D_{N}\Psi(D_{N_{P}}(N_{Q}\oplus\cdot)(Y)) &= D_{P}(Q\otimes\cdot)(V).
\end{align}
Altogether, \ref{pull}, \ref{bi}, \ref{X}, \ref{Y} imply,
\begin{multline*} 
    g_{P}(W,V) =\tilde{g}_{N_{P}}(X,Y) = \tilde{g}_{N_{Q}\oplus N_{P}}(D_{N_{P}}(N_{Q}\oplus\cdot)(X),D_{N_{P}}(N_{Q}\oplus\cdot)(Y))\\
    =g_{Q\otimes P}(D_{P}(Q\otimes\cdot)(W),D_{P}(Q\otimes\cdot)(V)),
\end{multline*}
as desired.

\subsection{Parallel transport} \label{subsection 5.2}

\noindent Let $N,K\in\mathcal{L}^{n\times p}_{*}$. Let $\tau_{N,K}$ denote the parallel transporter from $N$ to $K$ in $\mathcal{L}^{n\times p}_{*}$. Then verbatim from {\bf Lemma} \ref{lem:lin}, we have:
\begin{equation*}
    \tau_{N,K}: T_{N}\mathcal{L}^{n\times p}_{*}\ni X\mapsto (D_{N}((K\oplus N^{-1})\oplus\cdot)X\in T_{K}\mathcal{L}^{n\times p}_{*}.
\end{equation*}
Let $Y=(D_{N}((K\oplus N^{-1})\oplus\cdot)X$, then from \ref{dis}, $Y_{ii}= K_{ii}N_{ii}^{-1}X_{ii}, i=1,\cdots,p$ and $Y_{ij}=X_{ij}, i>j$.\\

\noindent
Similarly, let $P,Q\in S(n,p)^{*}$ and $W\in T_{P}S(n,p)^{*}$. Let $\tau_{P,Q}$ be the parallel transporter from $P$ to $Q$ in $S(n,p)^{*}$. Then, 
\begin{equation*}
    \tau_{P,Q}: T_{P}S(n,p)^{*}\ni W\mapsto (D_{P}((Q\otimes P^{-1})\otimes\cdot)W\in T_{Q}S(n,p)^{*}.
\end{equation*}
Now $((Q\otimes P^{-1})\otimes\cdot)$ expresses the left translation by $Q\otimes P^{-1}$ in $S(n,p)^{*}$. From \ref{obs}, we have a simpler expression for $\tau_{P,Q}(W)$ in terms of $n\times p$ matrices:
\begin{equation} \label{par} 
    (D_{P}((Q\otimes P^{-1})\otimes\cdot)W= (D_{N_{P}}((N_{Q}\oplus N_{P}^{-1})\oplus\cdot)(X))N_{Q}^{T}+N_{Q}(D_{N_{P}}((N_{Q}\oplus N_{P}^{-1})\oplus\cdot)(X))^{T},
\end{equation}
where $N_{P}N_{P}^{T}=P, N_{Q}N_{Q}^{T}=Q$, and $X=(D_{N_{P}}\Psi)^{-1}(W)=[N^{au}_{P}[(N^{au})^{-1}_{P}W(N^{au})^{-T}_{P}]_{1/2}]^{dr}$, and derivatives of translations in $\mathcal{L}^{n\times p}_{*}$ have an easy form as in \ref{dis}.

\subsection{A sample calculation} \label{subsection 5.3}

\noindent
Let $P,Q\in S(n,p)^{*}$ and $W\in T_{P}S(n,p)^{*}$ be as above. Then $\tau_{P,Q}(W)$, the translation of the tangent $W$ at $P$ to the tangent $\tau_{P,Q}(W)$ at $Q$ is given by \ref{par}. Now,
\begin{equation} \label{cal} 
    \begin{split}
        D_{N_{P}}((N_{Q}\oplus N_{P}^{-1})\oplus\cdot)(X))_{ii} &=(N_{Q})_{ii}(N_{P})^{-1}_{ii}X_{ii},\,\,i=1,\cdots,p\\
        (D_{N_{P}}((N_{Q}\oplus N_{P}^{-1})\oplus\cdot)(X))_{ij} &=X_{ij},\,\,i>j
    \end{split}
\end{equation}
as in \ref{dis}. Suppose
\begin{equation*}
    P=\begin{pmatrix}1&2&3\\2&5&8\\3&8&13\end{pmatrix},\,\, Q=\begin{pmatrix}1&2&-1\\2&5&-3\\-1&-3&2\end{pmatrix},\,\, \text{ and } \,\,W=\begin{pmatrix}2&3&4\\3&4&6\\4&6&10\end{pmatrix}.
\end{equation*}
Then,
\begin{equation*}
    N_{P}=\begin{pmatrix}1&0\\2&1\\3&2\end{pmatrix},\,\,N_{Q}=\begin{pmatrix}1&0\\2&1\\-1&-1\end{pmatrix},\,\, \text{ and } \,\,N^{au}_{P}=\begin{pmatrix}1&0&0\\2&1&0\\3&2&1\end{pmatrix}.
\end{equation*}
Then from \ref{indiff}, $X=(D_{N_{P}}\Psi)^{-1}(W)$ is,
\begin{equation*}
    X=[N^{au}_{P}[(N^{au}_{P})^{-1}W(N^{au}_{P})^{-T}]_{1/2}]^{dr}=\bigg[\begin{pmatrix}1&0&0\\2&1&0\\3&2&1\end{pmatrix}\cdot\begin{pmatrix}1&0&0\\-1&0&0\\0&1&0\end{pmatrix}\bigg]^{dr}=\begin{pmatrix}1&0\\1&0\\1&1\end{pmatrix}.
\end{equation*}
From \ref{cal}, the vector $K=D_{N_{P}}((N_{Q}\oplus N_{P}^{-1})\oplus\cdot)(X)$ is $K=\begin{pmatrix}1&0\\1&0\\1&1\end{pmatrix}$. Hence the tangent $\tau_{P,Q}(W)$ is,
\begin{equation} \label{parex} 
    \tau_{P,Q}(W)=KN_{Q}^{T}+N_{Q}K^{T}=\begin{pmatrix}2&3&0\\3&4&2\\0&2&-4\end{pmatrix}.
\end{equation}

\section{Reduced Cholesky factorization and LRC algorithm} \label{section 6}

\subsection{An algorithm for reduced Cholesky factorization} \label{subsection 6.1}

\noindent
A quick reduced Cholesky factorization for $M\in S(n,p)^{*}$ can be done as follows. Let $\Lambda\in\mathbb{R}^{p\times p}$ be the positive spectrum of $M$; ie, $\Lambda$ is the diagonal matrix whose diagonal entries are all the positive eigenvalues of $M$ listed decreasingly (the ordering of these eigenvalues will not affect the algorithm, actually). Let 
\begin{equation*}
    M = U\Lambda U^{T}=(U\Lambda^{1/2})(U\Lambda^{1/2})^{T}
\end{equation*}
be a thin spectrum decomposition of $M$. Perform the QR decomposition on $(U\Lambda^{1/2})^{T}=\Lambda^{1/2}U^{T}$ to have an orthogonal matrix $Q\in\mathbb{R}^{p\times p}$ and a mock upper triangular matrix $R\in\mathbb{R}^{p\times n}$ whose mock diagonal entries are positive ($Q$, $R$ will be unique, as argued below). Then
\begin{equation*} 
    M=(U\Lambda^{1/2})(U\Lambda^{1/2})^{T}= R^{T}Q^{T}QR = R^{T}R.
\end{equation*}
Set $N=R^{T}$. Then $N$ is the desired reduced Cholesky factor of $M$. This algorithm works as intended due to the following reason. Let $S(n,p)^{*}\ni M = ZZ^{T}$ where $Z=U\Lambda^{1/2}\in\mathbb{R}^{n\times p}_{*}$:
\begin{equation*} 
    M=ZZ^{T}=Z\,\,\bigg[ col(Z^{T})_1\,\,col(Z^{T})_2\,\cdots\,col(Z^{T})_{p}\,\,\cdots\,\bigg] = \bigg[col(M)_1\,\,col(M)_2\,\cdots\,col(M)_{p}\,\,\cdots\,\bigg].
\end{equation*}
If linear dependence appears among $col(Z^{T})_1,\cdots,col(Z^{T})_{p}$, the same is true for $col(M)_1,\cdots,col(M)_{p}$, which is a contradiction. Hence the QR decomposition for $Z^{T}$ yields:
\begin{equation} \label{QR} 
    \Lambda^{1/2}U^{T}=Z^{T}=QR = Q[R_1\,\,R_2],
\end{equation}
where $QR_1$ is precisely the unique QR decomposition of the full-rank square matrix \cite{golub2013matrix}
\begin{equation*}
    \bigg[ col(Z^{T})_1\,\,col(Z^{T})_2\,\cdots\,col(Z^{T})_{p}\bigg]\in\mathbb{R}^{p\times p},
\end{equation*}
with $(R_1)_{ii}>0, i=1,\cdots,p$, and $R_2\in\mathbb{R}^{p\times(n-p)}$ be such that $QR_2$ forms the remaining columns of $Z^{T}$.\\

\noindent
Note that through this algorithm, one can easily get the positive spectrum $\Lambda$ of $M$ back through $N$, as $U\Lambda^{1/2}=R^{T}Q^{T}=NQ^{T}$ and hence,
\begin{equation*}
    N^{T}N=(U\Lambda^{1/2}Q)^{T}(U\Lambda^{1/2}Q)=Q^{T}\Lambda Q\in P(p).
\end{equation*}

\subsection{LRC algorithm} \label{subsection 6.2}

We introduce our Low Rank Cholesky (LRC) algorithm in Algorithm 1.

\begin{algorithm}[h]
\caption{Mean via Low Rank Cholesky (LRC) algorithm}
\label{alg:LRC}
\textbf{Input :} $n\times n$ matrices $\{ B_i \}_{i = 1}^m$ with rank $K$, $n>K$\\
\textbf{Step 1:} Compute the reduced Cholesky factor $L_i\in\mathbb{R}^{n \times K}$ of $B_i$, as described in subsection \ref{subsection 6.1} and subsection \ref{subsection 9.3}. \\
\textbf{Step 2:} Compute the mean $\mathcal{F}(L_1,\cdots,L_m)$ as in \ref{mean2}. \\
\textbf{Output :} $B=\mathcal{F}(L_1,\cdots,L_m)\mathcal{F}(L_1,\cdots,L_m)^{T}\in \mathbb{R}^{n \times n}$
\end{algorithm}

\noindent
We will use this algorithm in section \ref{section 7} to estimate eigenspaces of an underlying low-rank covariance matrix from a set of random positive semidefinite matrices.  

\section{Applications} \label{section 7}

\noindent
Let $\Sigma$ denote the underlying true low-rank covariance matrix from a set of distributed random samples. Through eigendecomposition, we obtain $\Sigma = V \Lambda V^T$, where $\Lambda=\text{diag}(\lambda_1, \ldots, \lambda_n)$ with $\lambda_1 \geq \ldots \geq \lambda_n$, and $V = (v_1, \cdots, v_n)$ stores corresponding orthonormal eigenvectors. Let $V_K = (v_1, \cdots, v_K)$. We demonstrate the usage of our LRC algorithm in estimating $\text{Col}(V_K)$, the linear space spanned by the top $K$ eigenvectors of $\Sigma$ from the distributed data. This has many applications in machine learning literature \cite{fan2019distributed, liang2014improved, qu2002principal}. We investigate two scenarios with synthetic data, when the input are: 1) $m$ noisy $S(n, K)^*$ matrices with the underlying true matrix being $\Sigma$ in Section \ref{subsection 7.1}, and 2) $m$ sets of $N$ i.i.d. random samples $\{ X_i\}_{i=1}^N \subset \mathbb{R}^{n}$ with mean zero and covariance matrix $\Sigma$ in Section \ref{subsection 7.2}. To measure the estimation accuracy, we use $\rho(\hat{V}_K, V_K) = \|\hat{V}_K \hat{V}_K^T - V_K V_K^T \|$, which is also known as the $\sin\Theta$ distance between $\hat{V}_k$ and $V_k$. 

\subsection{Estimate principal eigenspace from noisy matrices} \label{subsection 7.1}

\noindent In this section, we generate the underlying $\Sigma \in S(n, K)^*$ with $K = 3$ as follows. We first create $W=(v_1, v_2, v_3) \in \mathbb{R}^{n \times 3}$ with the first four elements of $(v_1, v_2, v_3)$ being
\begin{equation*} 
    (1/2, 1/2, 1/2, 1/2)^{T}, (-1/2, -1/2, 1/2, 1/2)^{T}, (1/2, -1/2, 1/2, -1/2)^{T}, 
\end{equation*}
respectively and the rest being zero. We then obtain the matrix $\Sigma$ through $\Sigma = W  \Lambda_3 W^T$, where $\Lambda_3 =$ diag$( \lambda_1, \lambda_2, \lambda_3)$. To introduce randomness to the eigenvalues, we generate  $\lambda_1 \sim N(10, 1^2), \lambda_2 \sim N(5, 1^2)$,  and $ \lambda_3 \sim N(2.5, 1^2)$. 
%Since $\Sigma$ is rank 3, its reduced Cholesky decomposition gives, $\Sigma = \mathbf{N}_3\mathbf{N}_3^T$, where  $\mathbf{N}_3 \in \mathbb{R}^{n\times 3}$ is a mock lower triangular matrix.
To generate the $i$th noisy matrix data, we perturb $\Sigma$ through 
\[ \Sigma^{(i)} = \Sigma +  (E^{(i)} + E^{(i)T})/2,
\] where $E^{(i)}_{(l, j)} \sim N(0, \sigma^2_\epsilon)$. We consider twelve settings, with $n\in\{50,200\}$, corresponding respectively to low and high ambient dimensions, $\sigma_\epsilon^2 \in \{0.1^2, 0.3^2, 0.5^2 \}$ corresponding to different noise levels, and  $m \in \{ 20,100 \}$ matrices, corresponding respectively to small and large sample sizes.
Besides our LRC method, we also utilize distributed PCA (dPCA) proposed by \cite{fan2019distributed} and an ad-hoc Eigenvector averaging (EigenV-ave) method to estimate $\text{Col}(V_K)$. Briefly, for the EigenV-ave method, we first calculate the top $K$ eigenvectors $\{ \hat{V}_K^{(i)}\}_{i = 1}^m$ of each noisy matrix $\Sigma^{(i)}$, then compute the average $\tilde{V}_K = (1/m)\sum_{i= 1}^m \hat{V}_K^{(i)}$ and $\text{Col}(\tilde{V}_K)$.\\

\begin{figure}[h]
    \centering
    \includegraphics[scale=0.9]{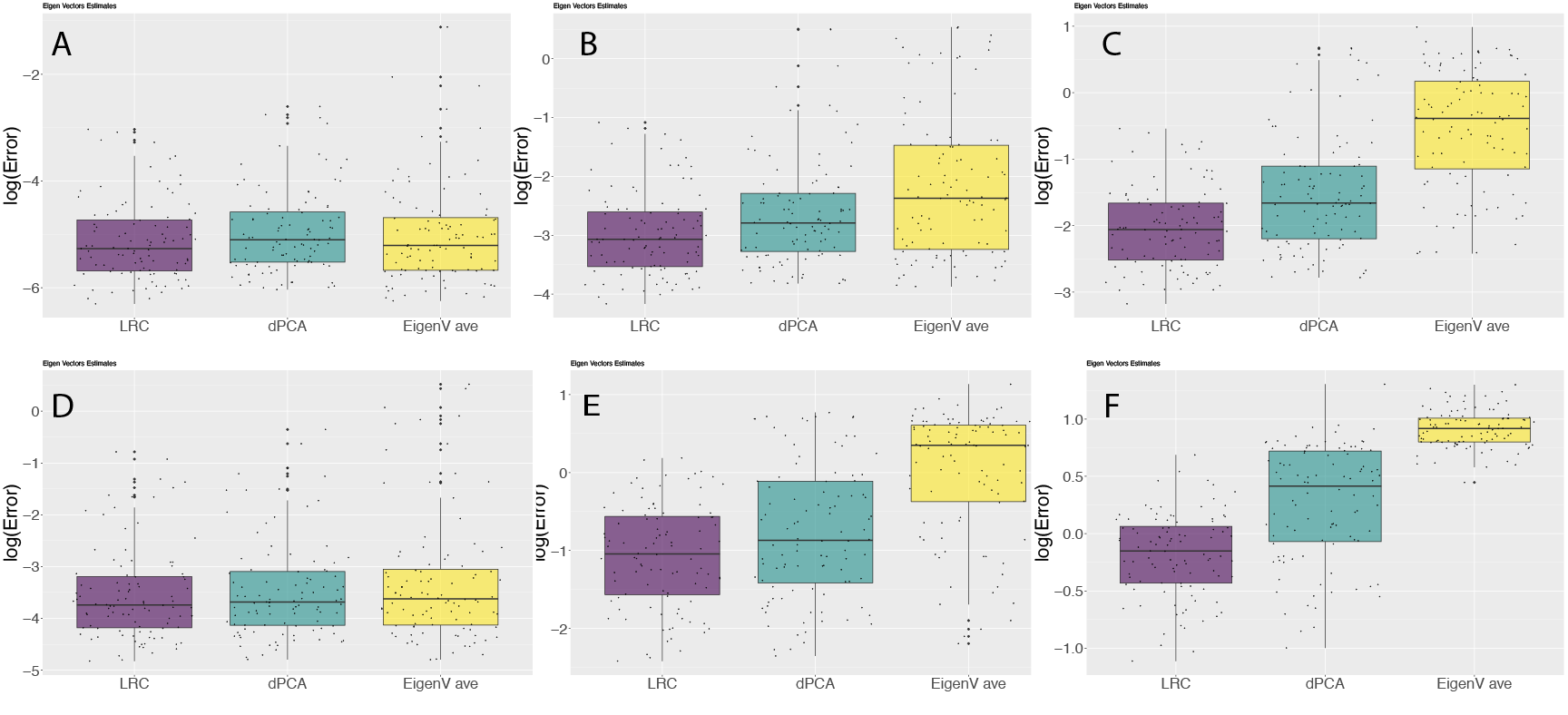}
    \caption{Comparison between LRC, dPCA, and EigenV-ave using perturbed random matrices: (A) $\sigma_\epsilon^2 = 0.1$ and $n = 50$; (B) $\sigma_\epsilon^2 = 0.3$ and $n = 50$; (C) $\sigma_\epsilon^2 = 0.5$ and $n = 50$; (D) $\sigma_\epsilon^2 = 0.1$ and $n = 200$; (E) $\sigma_\epsilon^2 = 0.3$ and $n = 200$; (F) $\sigma_\epsilon^2 = 0.5$ and $n = 200$}
    \label{fig:Fig3}
\end{figure}

\noindent
The results of the simulations based on 100 replicate trails are reported in Figure \ref{fig:Fig3} and Table \ref{tab:randommatrix}. When the noise level is small, all three methods have similar performance as shown in Figures \ref{fig:Fig3}A, \ref{fig:Fig3}D, since the noise to signal ratio is very low. LRC starts to gain advantage over dPCA and EigenV ave when $\sigma_\epsilon^2$ increases as in Figures \ref{fig:Fig3}B, \ref{fig:Fig3}C and \ref{fig:Fig3}E, \ref{fig:Fig3}F. The advantage persists under different dimensions ($n \in \{ 50, 200\}$) and different sample size ($m \in \{20, 100 \}$). In sum, our method dominates the other two methods for all twelve settings, especially when the noise level is large.\\

\subsection{Estimate principal eigenspace from random vectors} \label{subsection 7.2}

\noindent
Positive semidefinite matrices occur in numerous applications, including covariance matrices for multivariate measurements. Unlike the previous subsection where these matrices are given as inputs, we will calculate them from i.i.d. random vectors in this section. We assume the data are stored in $m$ different machines, each of which contains either $100$ or $50$ samples following $N(0, \Sigma_1)$, for some $\Sigma_1$. In fact, this is the distributed PCA setting \cite{fan2019distributed}. This setting is commonly used to handle scenarios in which large datasets are scattered across distant places, causing concerns for high communication costs and data security (for example, health records from different hospitals with strong privacy regulations). In these scenarios, summary statistics such as eigenvectors or Cholesky factors are preferred to the original data for aggregating for the final estimate. We first generate $\Sigma$ as in subsection \ref{subsection 7.1} and obtain the covariance matrix $\Sigma_1 = \Sigma + \sigma_\epsilon^2 I$, where $I$ is the $n \times n$ identity matrix and $\sigma_\epsilon^2 I$ corresponds to the covariance matrix of the measurement errors, which are independent of the signals, whose covariance matrix is $\Sigma$. We again consider two different dimensions $n \in \{50, 200 \}$, two different numbers of sites, $m \in \{ 20, 100 \}$, two different number of random vectors in each site, $l \in \{ 50, 100\}$, and three different noise level, $\sigma_\epsilon^2 \in \{0.1^2,0.3^2, 0.5^2\}$. Besides the LRC, dPCA, and EigenV-ave methods, we also consider the full sample PCA (fPCA) which pools all the samples together and hence is the oracle method.\\ 

\begin{figure}[h]
    \centering
    \includegraphics[scale=0.9]{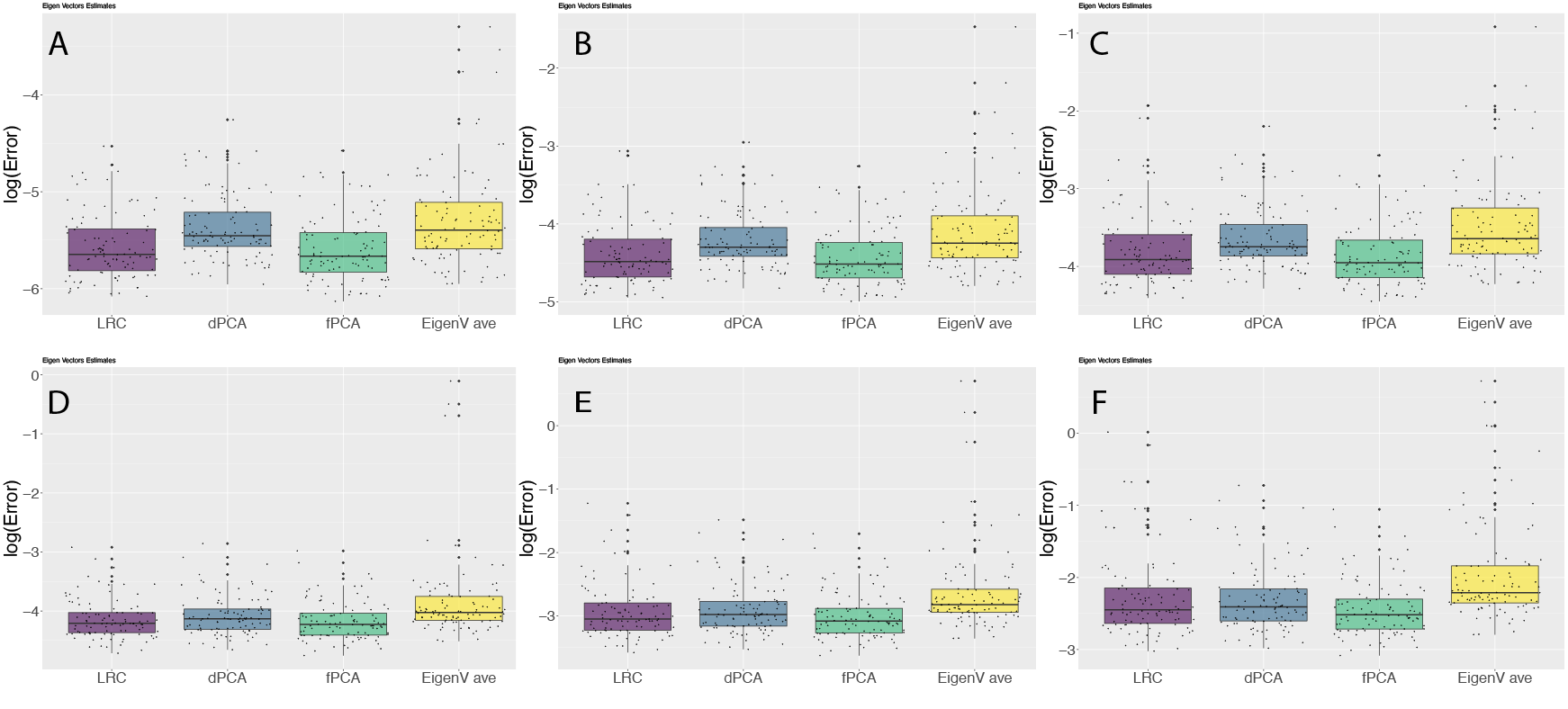}
    \caption{Comparison between LRC, dPCA, fPCA and EigenV-ave using random vectors: (A)  $\sigma_\epsilon^2 = 0.1$ and $n = 50$; (B)  $\sigma_\epsilon^2 = 0.3$ and $n = 50$; (C)  $\sigma_\epsilon^2 = 0.5$ and $n = 50$; (D)  $\sigma_\epsilon^2 = 0.1$ and $n = 200$; (E)  $\sigma_\epsilon^2 = 0.3$ and $n = 200$; (F)  $\sigma_\epsilon^2 = 0.5$ and $n = 200$}
    \label{fig:dPCA}
\end{figure}

\noindent
From Figures \ref{fig:dPCA}A-\ref{fig:dPCA}F and Table \ref{tab:randomvec}, we can see that fPCA has the best performance as expected. As shown in Figures \ref{fig:dPCA}A-\ref{fig:dPCA}C and Table \ref{tab:randomvec}, the LRC method has better performance than the dPCA and EigenV ave in all dimensional settings across different noise levels. When dimension increases, LRC has similar but slightly better estimation accuracy than dPCA and EigenV ave as shown in Figures \ref{fig:dPCA}D-\ref{fig:dPCA}F. 

\section{Conclusion}\label{section 8}

\noindent
In this paper, we introduced and investigated the special submanifold $S(n,p)^{*}$ of the manifold of positive semidefinite matrices, $S(n,p)$. We argued that it's a good substitute for $S(n,p)$ as it's a manifold of the same dimension $np-p(p-1)/2$ as $S(n,p)$ and that it's dense in $S(n,p)$ in the Frobenius sense. We endowed $S(n,p)^{*}$ with a Riemannian structure that would allow closed forms of various applied tools: endpoint geodesics, exponential and logarithmic maps, parallel transports and Fr\'echet means. We provided an algorithm to compute the reduced Cholesky factor for a matrix in $S(n,p)^{*}$. Finally, for applications, we demonstrated that our method for distributed subspace learning - via using the introduced LRC algorithm - surpassed other current methods. \\

\vspace{-20pt}
\begin{table}[htb]
\caption{ \label{tab:randommatrix} Summary statistics reporting performance of the LRC, dPCA and EigenV ave methods inferring the true eigenspace using random matrices. In each cell, the number is the mean of the estimation error based on 100 replicates, and the number in the parentheses is the corresponding standard deviation.}
\textcolor{blue}{
\begin{center}
  \scalebox{0.85}{
  \begin{tabular}{c c c c c c }
  \hline
 \thead{Dimension $(n)$}  & \thead{Sample Size $(m)$} & \thead{Noise level $(\sigma^2_\epsilon)$}  & \thead{LRC} & \thead{dPCA} & \thead{EigenV ave} \\
  \hline
  \multirow{6}{*}{$n = 50$}
  &\multirow{3}{*}{20}
    &  $0.1^2$    &  \bf{0.008 (0.009)} &  0.010 (0.013) & 0.014 (0.037)      \\
       && $0.3^2$     & \bf{0.068 (0.064)} & 0.132 (0.252) &  0.267 (0.408) \\
    && $0.5^2$     &   \bf{0.157 (0.112)} & 0.350 (0.426) & 0.809 (0.595) \\
    \cline{2-6}
&\multirow{3}{*}{100}
    & $0.1^2$    &  \bf{0.002 (0.002)} &  0.002 (0.003)   &  0.002 (0.005\\
    && $0.3^2$       &  \bf{0.014 (0.013)} &  0.022 (0.030)      &  0.080 (0.156) \\
    && $0.5^2$    &  \bf{0.034 (0.027)} &  0.073 (0.110) &  0.345 (0.398)   \\
  \hline
   \multirow{6}{*}{$n = 200$}
  &\multirow{3}{*}{20}
    & $0.1^2$      &  \bf{0.051 (0.076)} &  0.059 (0.103)   &  0.117 (0.287)      \\
       && $0.3^2$     &  \bf{0.422 (0.268)} &  0.662 (0.582) &  1.278 (0.691) \\
    && $0.5^2$     &  \bf{0.888 (0.315)} &  1.464 (0.649) &  2.52 (0.418) \\
    \cline{2-6}
&\multirow{3}{*}{100}
    & $0.1^2$      &  \bf{0.011 (0.019)} &  0.013 (0.027)  &  0.037 (0.112) \\
    && $0.3^2$      & \bf{0.11 (0.09)}  &  0.22 (0.32) & 0.73 (0.63)  \\
    && $0.5^2$    & \bf{0.44 (0.23)}  & 1.09 (0.66) & 1.96 (0.39)  \\
  \hline
  \end{tabular}}
  \end{center}}
\end{table}

\begin{table}[htb]
\caption{\label{tab:randomvec}  
Summary statistics reporting performance of the LRC, dPCA, fPCA and EigenV ave methods inferring the true eigenspace using random vectors with shared covariance matrix from $m$ different sites. In each cell, the number is the mean of the estimation error based on 100 replicates, and the number in the parentheses is the corresponding standard deviation.}
\textcolor{blue}{
\begin{center}
  \scalebox{0.75}{
  \begin{tabular}{c c c c c c c c}
  \hline
  \thead{Dimension \\ ($n$) } & \thead{ Number of sites \\ $(m)$}  & \thead{Noise level \\($\sigma^2_\epsilon$)} &  \thead{Sample size \\ within a site  ($l$) } & \thead{LRC} & \thead{dPCA} & \thead{fPCA} & \thead{EigenV ave} \\
  \hline
  \multirow{12}{*}{$n = 50$}
  &\multirow{6}{*}{20}
    &  $0.1^2$  & 50 &  0.008 (0.003) &  0.010 (0.004) & 0.008 (0.003) & 0.064 (0.260)    \\
       && $0.1^2$ & 100   & 0.004 (0.002) & 0.005 (0.002) &  0.004 (0.002)& 0.006 (0.005)  \\
     && $0.3^2$ & 50   & 0.028 (0.014) & 0.033 (0.015) &  0.025 (0.011)& 0.115 (0.298)  \\
       && $0.3^2$ & 100  & 0.014 (0.007) & 0.016 (0.007) &  0.013 (0.006)& 0.022 (0.026) \\
    && $0.5^2$   & 50 &   0.053 (0.035) & 0.059 (0.029) & 0.044 (0.021) & 0.182 (0.339)\\
    && $0.5^2$   & 100 &  0.026 (0.019) & 0.029 (0.015) & 0.023 (0.011) & 0.040 (0.047) \\
    \cline{2-8}
&\multirow{6}{*}{100}
    &  $0.1^2$  & 50 &  0.002 (0.001) &  0.002 (0.001) & 0.002 (0.001) & 0.006 (0.025)  \\
    && $0.1^2$ & 100 & 0.001 (0.064) & 0.001 (0.001) &  0.001 (0.001)& 0.001 (0.002)   \\
    && $0.3^2$ & 50 & 0.006 (0.003) & 0.006 (0.003) &  0.005 (0.002) & 0.020 (0.062) \\
    && $0.3^2$ & 100 & 0.003 (0.001) & 0.003 (0.001) & 0.003 (0.0.001)& 0.005 (0.009)  \\
    && $0.5^2$ &.50 &  0.011 (0.008) & 0.012 (0.006) & 0.009 (0.004) & 0.057 (0.224) \\
    && $0.5^2$ & 100 & 0.005 (0.004) & 0.006 (0.003) & 0.005 (0.002) & 0.010 (0.022) \\
  \hline
   \multirow{12}{*}{$n = 200$}
  &\multirow{6}{*}{20}
    & $0.1^2$     & 50 &  0.035 (0.014) & 0.037 (0.103) & 0.032 (0.012) & 0.092 (0.203)  \\
       && $0.2^2$ & 100&  0.017 (0.007) & 0.018 (0.007) & 0.016 (0.007) & 0.040 (0.116) \\
       && $0.3^2$ & 50 &  0.101 (0.093) & 0.103 (0.058) & 0.092 (0.043) & 0.158 (0.317) \\
       && $0.3^2$ & 100&  0.047 (0.041) & 0.051 (0.031) & 0.046 (0.024) & 0.059 (0.239) \\
       && $0.5^2$ & 50 &   0.193 (0.219) & 0.182 (0.114) & 0.160 (0.079)& 0.326 (0.410)\\
       && $0.5^2$ & 100&   0.086 (0.144) & 0.090 (0.070) & 0.081 (0.047)& 0.109 (0.292)\\
    \cline{2-8}
&\multirow{6}{*}{100}
    & $0.1^2$     & 50 & 0.007 (0.003) &  0.008 (0.003)  &  0.007 (0.003) & 0.015 (0.026) \\
    && $0.2^2$    & 100& 0.003 (0.001) &  0.004 (0.002) & 0.003 (0.001)  & 0.025 (0.200) \\
    && $0.3^2$    & 50 & 0.020  (0.027)  &  0.021 (0.014) & 0.018 (0.0.011)  & 0.031 (0.115) \\
    && $0.3^2$    & 100& 0.010 (0.009)  &  0.010 (0.007) & 0.010 (0.005)  & 0.012 (0.164)\\
    && $0.5^2$    & 50 & 0.041 (0.076)  & 0.037 (0.032) & 0.032 (0.021) & 0.069 (0.222) \\
    && $0.5^2$    & 100& 0.018 (0.043)  &  0.019 (0.016) & 0.017 (0.010)  & 0.023 (0.220) \\
  \hline
  \end{tabular}}
  \end{center}}
\end{table}

\section{Appendix} \label{section 9}

\subsection{Cholesky factorization for positive semidefinite matrices}

\noindent
It's a known fact that for every $M\in S(n,p)$ there exists a unique lower triangular $L$ such that
\begin{equation*}
    M=LL^{T}
\end{equation*}
and that there are exactly $p$ positive diagonal entries and $n-p$ zero diagonal ones, which then makes $rank(L)=p=rank(M)$, as expected. Whenever a diagonal entry $L_{ii}=0$, its corresponding column $col(L)_{i}$ is also set to be a zero column. This factor $L$ of $M$ is called its Cholesky factor \cite{gentle2012numerical}. In this first subsection, we show that the positions of the zero diagonal entries (and consequently the zero columns) in $L$ can be easily deduced through the algorithm itself and are dictated by the algebraic relationship between the columns of $M$.

\subsubsection{Where dependence is spotted} \label{subsubsection 9.1.1}
\noindent
Let $M=LL^{T}$ be the Cholesky decomposition for $M\in S(n,p)$ with $L$ being the unique lower triangular Cholesky factor of $M$. We claim that the positions of the $n-p$ zero diagonals (and hence $n-p$ zero columns) in $L$ can be predicted as follows. Going from left to right, let $col(M)_{i}$ be the first column of $M$ that depends on the previous columns, $col(M)_1,\cdots,col(M)_{i-1}$. Then $L_{ii}=0$ (or equivalently $col(L)_{i}$ is a zero column \cite{gentle2012numerical}). As we do not see this discussed in the literature, we give a proof for this claim here.\\

\noindent
Let's start with some experimental cases first. Let's say $n=3,p=2$. It's easy to see from the Cholesky decomposition that if $col(M)_{i}$ is a zero column, for any $i$, then so is the $i$th column of its Cholesky factor $L$. Now suppose that $M\in S(3,2)$ such that $col(M)_2=t\cdot col(M)_1$ for some $t\not=0$, then
\begin{equation*}
    \begin{cases}M_{12}=tM_{11}=M_{21}\\ M_{22}=tM_{12}=M_{22}\\ M_{23}=tM_{13}=M_{32}\end{cases},
\end{equation*}
which then deduces that $M_{22}=t^2M_{11}$. By the Cholesky construction of $L$ \cite{gentle2012numerical}, one has that,
\begin{equation*}
    L_{22}^2= M_{22}-L_{21}^2=M_{22}-\frac{M_{21}^2}{M_{11}}=M_{22}-\frac{t^2M_{11}^2}{M_{11}}=M_{22}-t^2M_{11}=0.
\end{equation*}
This turns the whole second column of $L$ into a zero column. Now suppose that $col(M)_3=t\cdot col(M)_1+s\cdot col(M)_2$ for some $t^2+s^2>0$. Since $L$ is lower triangular, the third column of $L$ is
\begin{equation*}
    col(L)_3=\begin{pmatrix}0\\ 0\\ c\end{pmatrix}
\end{equation*}
for some scalar $c$. Let $\begin{pmatrix}a\\ b\\ c\end{pmatrix}$ denote the third column of $L^{T}$. From $LL^{T}=M$, one has
\begin{equation} \label{as}
    a\cdot col(L)_1+b\cdot col(L)_2+c\cdot col(L)_3=col(M)_3=t\cdot col(M)_1+ s\cdot col(M)_2.
\end{equation}
By way of construction, $col(M)_1, col(M)_2$ only depends on $col(L)_1,col(L)_2$, hence the last sum in \ref{as} can be rewritten as a linear combination of $col(L)_1,col(L_2)$: $u\cdot col(L)_1+v\cdot col(L)_2$, for some $u,v\in\mathbb{R}$. This means, from \ref{as} and the fact that $t^2+s^2>0$, $col(L)_3$ is a linear combination of $col(L)_1, col(L)_2$. As $L$ is a lower triangular matrix, this can only happen if $c=0$, hence the third column of $L$ is a zero column.\\

\noindent This last case above gives a hint on how to expand this proof further for the general case, as \ref{as} remains true if $i$ ($i>1$) substitutes for $3$, and $1,\cdots,i-1$ substitute for $1,2$. Suppose $col(M)_{i}$ is a linear combination of $col(M)_1,\cdots,col(M)_{i-1}$. Then from $LL^{T}=M$, one still arrives at the same conclusion as before:
\begin{equation} \label{as1} 
    \sum_{j=1}^{i} t_{j}\cdot col(L)_{j} = \sum_{j=1}^{i-1} a_{j}\cdot col(L)_{j}\Rightarrow t_{i}\cdot col(L)_{i} = \sum_{j=1}^{i-1} u_{j}\cdot col(L)_{j}
\end{equation}
for some $t_{j},a_{j},u_{j}\in\mathbb{R}$, where $t_{i}$ is also precisely the $i$th entry (top down) in $col(L)_{i}$. If $t_{i}\not= 0$ then \ref{as1} implies that either $col(L)_{i}=\vec{0}$ or that $col(L)_{i}$ is a nonzero vector and a linear combination of the previous columns $col(L)_1,\cdots,col(L)_{i-1}$. Either way it leads to a contradiction. Hence $t_{i}=0$ and the whole column $col(L)_{i}$ is a zero column.\\

\noindent For the other way around, it's easy to see that, via an induction proof, if the Cholesky factor $L$ of $M$ is such that $col(L)_{i}=\vec{0}$ then $col(M)_{i}$ is a linear combination of $col(M)_1,\cdots, col(M)_{i-1}$.

\subsection{Continuity of the mapping $M\mapsto L$} \label{subsection 9.2}

\noindent Let $\mathcal{L}$ be the set of all Cholesky factors of matrices in $S(n,p)$. In other words, $\mathcal{L}$ consists of all $n\times n$ lower triangular matrices with exactly $p$ positive diagonal entries and $n-p$ zero columns. As to be consistent with the terminology in \cite{lin2019riemannian}, we call this set, the Cholesky space (of $S(n,p)$). (Note that in \cite{lin2019riemannian}, $n=p$, and the Cholesky space there is of $P(n)$, denoted $\mathcal{L}_{+}$, and consists of $n\times n$ lower triangular matrices with all positive diagonal entries.)\\
Let $\Phi:\mathcal{L}\to S(n,p)$ where,
\begin{equation*}
    \Phi(L)=LL^{T} = M
\end{equation*}
denote the Cholesky map, and let $\Phi^{-1}: S(n,p)\to\mathcal{L}$ where,
\begin{equation*}
    \Phi^{-1}(M) = L
\end{equation*}
with $LL^{T}=M$. The existence and uniqueness of the Cholesky decomposition make $\Phi^{-1}$ a well-defined map. However, to be clear, the map $\Phi^{-1}$ fails to be continuous when $M\in S(n,p)$, unlike the case when $M\in P(n)$ \cite{lin2019riemannian}. This is because the automatic setting of a column in $L$ to a zero column once its corresponding diagonal entry works out to be zero \cite{gentle2012numerical}. On the other hand, continuity is preserved if one has determined in advance which are the first $p$ linearly independent columns of $M$; as it's shown in sub-subsection \ref{subsubsection 9.1.1}, the corresponding columns in its Cholesky factor $L$ will be the only nonzero columns of $L$. In particular, if $M$ has its first $p$ columns linearly independent - i.e. $M\in S(n,p)^{*}$ - then $L$ looks like
\begin{equation} \label{ex1} 
    L=\begin{pmatrix}a_{11}&0&\dots&0&0&\dots&0\\ a_{21}&a_{22}&\dots&0&0&\dots&0\\ \vdots&\vdots&\ddots&\vdots&\vdots&\dots&\vdots\\ a_{p1}&a_{p2}&\dots&a_{pp}&0&\dots&0\\ \vdots&\vdots&\dots&\vdots&\vdots&\dots&\vdots\\ a_{n1}&a_{n2}&\dots&a_{np}&0&\dots&0\end{pmatrix},
\end{equation}
where $a_{ii}>0$, $i=1,\cdots,p$. Then from the fact that $LL^{T}=M$, the entries of $M$ can be solved continuously in terms of the entries of $L$ and vice versa. Let us illustrate this with $n=3,p=2$:
\begin{equation} \label{const}
    LL^{T}=\begin{pmatrix}a&0&0\\b&d&0\\c&e&0\end{pmatrix}\begin{pmatrix}a&b&c\\ 0&d&e\\0&0&0\end{pmatrix}=\begin{pmatrix}a^2&ab&ac\\ ab&b^2+d^2&bc+de\\ ac&bc+de&c^2+e^2\end{pmatrix}=\begin{pmatrix}M_{11}&M_{12}&M_{13}\\M_{21}&M_{22}&M_{23}\\M_{31}&M_{32}&M_{33}\end{pmatrix} = M.
\end{equation}
This is precisely the Cholesky map $\Phi: L\mapsto LL^{T}=M$. One can see that the entries of $M$ are polynomial mappings of the entries of $L$. Conversely, the entries of $L$ are constructed in terms of the entries of $M$ as follows \cite{gentle2012numerical}:
\begin{align} 
    \nonumber a &= \sqrt{a^2} = \sqrt{M_{11}}\\ 
    \nonumber b &= (ab)/a= M_{21}/\sqrt{M_{11}}\\ 
    \label{map1} c &= (ac)/a = M_{31}/\sqrt{M_{11}}\\
    \nonumber d &= \sqrt{(b^2+d^2)-b^2} = \sqrt{M_{22}-M_{21}^2/M_{11}}\\ 
    \nonumber e &= (bc+de -bc)/d = \big(M_{32}- (M_{21}M_{31})/M_{11}\big)/\sqrt{M_{22}-M_{21}^2/M_{11}}.
\end{align}
The diagonal entries $a,d$ of $L$ in \ref{const} are both positive. One can see from \ref{map1} that entries of $L$ are expressed as continuous functions of entries of $M$ (summation, multiplication, division or taking roots). This statement is true in general for $M\in S(n,p)^{*}$, as the nonzero entries of $L$ can be solved inductively through the Cholesky program \cite{gentle2012numerical}:
\begin{equation} \label{chp} 
    \begin{cases}L_{ii}=\sqrt{M_{ii}-\sum_{k=1}^{i-1}L^2_{ik}},\,\,i=1,\cdots,p\\
    L_{ij}=\frac{1}{L_{jj}}(M_{ij}-\sum_{k=1}^{j-1}L_{ik}L_{jk}),\,\, j=1,\cdots,p, j<i.
    \end{cases}
\end{equation}
This means, $\Phi,\Phi^{-1}$ are both continuous (in fact, smooth) and one-to-one maps, if restricted to matrices in $S(n,p)^{*}$ and their Cholesky factors. Smoothness is measured using the usual Euclidean metric in $\mathbb{R}^{n^2}\cong\mathbb{R}^{n\times n}$.\\

\noindent
There is more to \ref{ex1}. Since $\Phi^{-1}|_{S(n,p)^{*}}$ serves an identification between $M\in S(n,p)^{*}$ and its Cholesky factor $L$, one can then use \ref{ex1} to give a local Euclidean coordinate chart for $M\in S(n,p)^{*}$. More precisely:
\begin{equation} \label{1to1} 
    S(n,p)^{*}\ni M\overset{\Phi^{-1}}{\longleftrightarrow} L = (a_{ij})_{1\leq j\leq p, j\leq i}.
\end{equation}
Clearly, the tuples $(a_{ij})_{1\leq j\leq p, j\leq i}$, with $a_{ii}>0, i=1,\cdots,p$ and other $a_{ij}\in\mathbb{R}$, carve out open neighborhoods in $\mathbb{R}^{np-p(p-1)/2}$. Hence the identification \ref{1to1} gives a proof that $S(n,p)^{*}$ is itself a smooth manifold of dimension $np-p(p-1)/2$. 

\subsection{A submanifold} \label{subsection 9.3}

\noindent
The reasoning in subsection \ref{subsection 9.2} doesn't allow one to extend \ref{1to1} to a consistent local coordinate chart for the whole $S(n,p)$, as $L$ looks differently for $M$ in a different class $C_{\{i_1,\cdots,i_{p}\}}$ (recall this notation from section \ref{section 2}). Hence the discussion doesn't allow one to conclude that $S(n,p)^{*}$ is a submanifold of $S(n,p)$. We will use another route. We first invoke the following fact. 

\begin{lemma} \label{lem:guillemin} \cite{guillemin2010differential}
If $X, Z$ are both (smooth) manifolds in some Euclidean space $\mathbb{R}^{N}$ and if $Z\subset X$, then $Z$ is a submanifold of $X$.
\end{lemma}

\noindent
It's a known fact that $S(n,p)$ is a $C^{\infty}$ (smooth) embedded submanifold of $\mathbb{R}^{n\times n}$ \cite{vandereycken2013riemannian}. Hence if we can show that $S(n,p)^{*}$ is a (smooth) manifold in $\mathbb{R}^{n\times n}$, then by {\bf Lemma} \ref{lem:guillemin}, $S(n,p)^{*}$ is a submanifold of $S(n,p)$. In order to show that $S(n,p)^{*}$ is a smooth manifold in $\mathbb{R}^{n\times n}$ of dimension $np-p(p-1)/2$, roughly speaking, we need to show that there exist a local smooth map $\Psi$, mapping from a neighborhood in $\mathbb{R}^{np-p(p-1)/2}$ to a neighborhood in $S(n,p)^{*}$, and a local smooth map $\Xi$, mapping from a neighborhood in $\mathbb{R}^{n\times n}$ to a neighborhood in $\mathbb{R}^{np-p(p-1)/2}$, such that, $\Psi$ and $\Xi|_{S(n,p)^{*}}$ are inverses of each other (\cite{guillemin2010differential}).\\

\noindent
Firstly, we note that if $M\in S(n,p)$ and $L$ is its Cholesky factor, then the zero columns of $L$ and correspondingly the zero rows of $L^{T}$ don't contribute to the calculation of $M=LL^{T}$. Let $N$ be the reduction of $L$ by throwing out the zero columns - call it the {\it reduced} Cholesky factor of $M$. Then $N\in\mathbb{R}^{n\times p}_{*}$ is a mock lower triangular matrix whose mock diagonal entries $N_{ii}\geq 0, i=1,\cdots,p$. For instance, if $n=4, p=2$ and $M\in C_{\{1,4\}}\subset S(4,2)$, then $L, N$ look like:
\begin{equation} \label{sampN} 
    L =\begin{pmatrix}a_{11}&0&0&0\\a_{21}&0&0&0\\a_{31}&0&0&0\\a_{41}&0&0&a_{44}\end{pmatrix}\rightarrow N=\begin{pmatrix}a_{11}&0\\a_{21}&0\\a_{31}&0\\a_{41}&a_{44}\end{pmatrix},
\end{equation}
where $a_{11},a_{44}>0$. Certainly, $LL^{T} = NN^{T}$. It follows from sub-subsection \ref{subsubsection 9.1.1} that if $M\in S(n,p)^{*}$ then $N_{ii}>0, i=1,\cdots,p$.\\

\noindent
Let $\mathcal{L}^{n\times p}$ denote the vector space of $n\times p$ mock lower triangular matrices and $\mathcal{L}^{n\times p}_{*}$ the set of all reduced Cholesky factors coming from $S(n,p)^{*}$ - we call the latter the {\it reduced} Cholesky space. Then $\mathcal{L}^{n\times p}_{*}\subset\mathcal{L}^{n\times p}$ as an open subset, and if $N\in\mathcal{L}^{n\times p}_{*}$ then $N_{ii}>0, i=1,\cdots,p$ and conversely. Moreover, $\mathcal{L}^{n\times p}\cong\mathbb{R}^{np-p(p-1)/2}$ as vector spaces, and $\mathcal{L}^{n\times p}_{*}\cong\mathcal{L}^{n\times p}\cong\mathbb{R}^{np-p(p-1)/2}$ as manifolds. Consider the two Euclidean spaces, $\mathcal{L}^{n\times p}\cong\mathbb{R}^{np-p(p-1)/2}$ and $\mathbb{R}^{n\times n}$, with their usual Euclidean metrics. Let $\mathcal{U}$ denote an open ball of $\mathcal{L}^{n\times p}$ such that if $N\in\mathcal{U}$ then $N_{ii}>0, i=1,\cdots,p$; ie, $\mathcal{U}\subset\mathcal{L}^{n\times p}_{*}$. Let 
\begin{equation*}
    \Psi:\mathbb{R}^{np-p(p-1)/2}\supset\mathcal{U}\to S(n,p)^{*}\subset\mathbb{R}^{n\times n}
\end{equation*} 
be such that $\Psi(N) = NN^{T}$. By the previous discussions, this map $\Psi$ is a well-defined smooth map. We want to construct a local $\Psi^{-1}$. For the sake of visualization, suppose $n=4, p=2$ - the following logic is independent of the choices of $n,p$. Then a typical image $\Psi(N)=NN^{T}$ looks like:
\begin{equation} \label{sampM} 
    \begin{split}
        \Psi(N) = NN^{T} = \begin{pmatrix}a&0\\b&e\\c&f\\d&g\end{pmatrix}\begin{pmatrix}a&0\\b&e\\c&f\\d&g\end{pmatrix}^{T} &= \begin{pmatrix}a^2&ab&ca&ad\\ab&b^2+e^2&bc+ef&bd+eg\\ac&bc+ef&c^2+f^2&cd+fg\\ad&bd+eg&cd+fg&d^2+g^2\end{pmatrix} \\
        &=\begin{pmatrix} M_{11}&M_{12}&M_{13}&M_{14}\\ M_{21}&M_{22}&M_{23}&M_{24}\\M_{31}&M_{32}&M_{33}&M_{34}\\M_{41}&M_{42}&M_{43}&M_{44}\end{pmatrix} = M.
    \end{split}
\end{equation}
The inversion $M\mapsto N$ of \ref{sampM} clearly follows from the inversion of the Cholesky map, with $L$ now being replaced by $N$. In the case of \ref{sampM}, the inversion looks like:
\begin{equation} \label{inv1} 
    N = \begin{pmatrix}\sqrt{M_{11}}&0\\M_{21}/\sqrt{M_{11}}&\sqrt{M_{22}-M_{21}^2/M_{11}}\\M_{31}/\sqrt{M_{11}}&(M_{32}-(M_{21}M_{31})/M_{11})/\sqrt{M_{22}-M_{21}^2/M_{11}}\\ M_{41}/\sqrt{M_{11}}&(M_{42}-(M_{21}M_{41})/M_{11})/\sqrt{M_{22}-M_{21}^2/M_{11}}\end{pmatrix}.
\end{equation}
The only complications for the inversion process are the necessary divisions in the Cholesky algorithm. However, because $N\in\mathcal{U}$ is such that $a,e>0$, each image $\Psi(N)=M$ satisfies
\begin{equation} \label{inv2} 
    M_{11}>0\text{ and } \sqrt{M_{22}-M_{21}^2/M_{11}}>0;
\end{equation}
i.e., the divisions in \ref{inv1} are defined. For each such $M=\Psi(N)$, let $B_2(M,\epsilon_{M})$ be an $l^2$ ball centered at $M$ with radius $\epsilon_{M}$ small enough so that the open conditions \ref{inv2} are satisfied everywhere in $B_2(M,\epsilon_{M})$. Let $\mathcal{V}=\bigcup_{M=\Psi(N); N\in\mathcal{U}} B_2(M,\epsilon_{M})$, which is an open set of $\mathbb{R}^{4\times 4}$. Define the smooth map $\Xi: \mathcal{V}\to\mathcal{L}^{4\times 2}$ as follows,
\begin{equation} \label{inv3}
    \begin{split}
        \Xi(M) &=\Xi(M_{11},M_{21},\cdots,M_{34},M_{44}) \\ 
        &=\begin{pmatrix} \sqrt{M_{11}}&0\\
        M_{21}/\sqrt{M_{11}}&\sqrt{M_{22}-M_{21}^2/M_{11}}\\M_{31}/\sqrt{M_{11}}&(M_{32}-(M_{21}M_{31})/M_{11})/\sqrt{M_{22}-M_{21}^2/M_{11}}\\M_{41}/\sqrt{M_{11}}&(M_{42}-(M_{21}M_{41})/M_{11})/\sqrt{M_{22}-M_{21}^2/M_{11}}\end{pmatrix}.
    \end{split}
\end{equation}
In other words, $\Xi(M)=\Xi(M')$ if $M_{i1}=(M')_{i1},i=1,\cdots,4, M_{i2}=(M')_{i2}, i=2,\cdots,4$, and $\Xi$ is a function of only $np-p(p-1)/2=7$ entries of $M$. By construction, $\Xi$ is well-defined on $\mathcal{V}$, and that from \ref{sampM}, \ref{inv3}
\begin{equation} \label{inv4} 
    \Psi\circ\Xi|_{\mathcal{V}\cap S(4,2)^{*}} = Id_{\mathcal{V}\cap S(4,2)^{*}} \,\,\text{ and }\,\, \Xi|_{\mathcal{V}\cap S(4,2)^{*}}\circ\Psi = Id_{\mathcal{U}}.
\end{equation}
This is what we want. As $\mathcal{V}\cap S(4,2)^{*}$ is an open set of $S(4,2)^{*}$, \ref{inv4} says that $\Psi$ is a local parametrization of $S(4,2)^{*}$ while $\Xi$ is its local coordinate chart \cite{guillemin2010differential}. Using the terminologies in \cite{guillemin2010differential}, \ref{inv4} says that $S(4,2)^{*}$ is a smooth manifold of dimension $np-p(p-1)/2=7$ in $\mathbb{R}^{4\times 4}$. \\

\noindent
While this concrete choice of $n=4,p=2$ was made to get the point (\ref{inv4}) across, it should be noted again that this logic carries out the same way for any $n,p$. The positive inequalities placed on the mock diagonal entries $N_{ii}, i=1,\cdots,p$ in \ref{inv2} are still strict inequalities ($N\in\mathcal{U}$); they are solved inductively in terms of elements in $M$ as follows \cite{gentle2012numerical}:
\begin{equation} \label{gen1} 
    0< N_{ii}=\sqrt{M_{ii}-\sum_{k=1}^{i-1}N_{ik}^2},\,\,i=1,\cdots,p.
\end{equation}
If $n,p$ are fixed, then $N_{ii}$ in \ref{gen1} can be eventually expressed as summations, multiplications, divisions or roots of entries in $M$, in a unique manner \cite{gentle2012numerical}. There are only $n\times n$ entries in $M$; hence there exists a small open region in $\mathbb{R}^{n\times n}$ centered at $M$ such that the strict inequalities \ref{gen1} hold. Note that not all the elements $M_{ij}$ are involved in the inversion process - only those needed to make \ref{gen1}. Hence the choice of $\mathcal{V}$ is possible and so is the general form of $\Xi$ in \ref{inv3}. We still have the conclusion that $S(n,p)^{*}$ is a smooth manifold of dimension $np-p(p-1)/2$ in $\mathbb{R}^{n\times n}$.

\vspace{5pt}

\begin{remark}
The conditions in \ref{inv2}, \ref{gen1} are technically conditions for entries of a symmetric matrix \cite{gentle2012numerical}. However, it's not a problem to extend these out for $n\times n$ matrices, as they are expressed in terms of the lower half entries.
\end{remark}

\subsubsection{Openness} \label{subsubsection 9.3.1}

\noindent
Pick $\tilde{M}\in S(n,p)^{*}$. As the conditions in \ref{gen1} (\ref{inv2} for an example) are open conditions, there exists a small open ball $\mathbb{R}^{n\times n}\supset\mathcal{V}\ni\tilde{M}$ such that \ref{gen1} (or \ref{inv2}) is satisfied for every $M\in\mathcal{V}$. For instance, if $n=4,p=2$, then $M\in\mathcal{V}$ looks like:
\begin{equation*}            
    M=\begin{pmatrix}M_{11}&M_{12}&M_{13}&M_{14}\\ M_{21}&M_{22}&M_{23}&M_{24}\\M_{31}&M_{32}&M_{33}&M_{34}\\M_{41}&M_{42}&M_{43}&M_{44}\end{pmatrix},
\end{equation*}
with,
\begin{equation*}
    M_{11}>0\text{ and } \sqrt{M_{22}-M_{21}^2/M_{11}}>0.
\end{equation*}
Now $\mathcal{V}\cap S(n,p)$ will consist of positive semidefinite matrices in $S(n,p)$ with conditions \ref{gen1} met (and if $n=4,p=2$, \ref{inv2}). Since these conditions are the only conditions required to solve for nonzero entries of $N\in\mathcal{L}^{n\times p}_{*}$ - as the mock diagonal entries $N_{ii}>0, i=1,\cdots,p$ - this means $\mathcal{V}\cap S(n,p)\subset S(n,p)^{*}$. Hence $\mathcal{V}\cap S(n,p)\subset S(n,p)^{*}$ is an open neighborhood of $\tilde{M}\in S(n,p)^{*}$ in $S(n,p)$, and thus $S(n,p)^{*}$ is open in $S(n,p)$.

\subsubsection{Dimensions}

\noindent
One can consider a similar approach to that in subsection \ref{subsection 9.2}, \ref{subsection 9.3} with other classes in $S(n,p)$. Suppose for convenience $p>3$ and consider $C_{\{1,2,4,5,\cdots,p+1\}}\subset S(n,p)$ (recall this notation from section \ref{section 2}. Then for $M\in C_{\{1,2,4,5,\cdots,p+1\}}$, its reduced Cholesky factor has the form:
\begin{equation} \label{ex2} 
    N=\begin{pmatrix}a_{11}&0&0&\dots&0\\ a_{21}&a_{22}&0&\dots&0\\a_{31}&a_{32}&0&\dots&0\\a_{41}&a_{42}&a_{43}&\dots&0\\ \vdots&\vdots&\vdots&\ddots&\vdots\\ a_{p+1\,1}&a_{p+1\,2}&a_{p+1\,3}&\dots&a_{p+1\,p}\\ \vdots&\vdots&\vdots&\dots&\vdots\\ a_{n1}&a_{n2}&a_{n3}&\dots&a_{n\,p}\end{pmatrix}.
\end{equation}
If one follows the strategy in subsection \ref{subsection 9.2} and derives a formulation similar to \ref{1to1} for $C_{\{1,2,4,5,\cdots,p+1\}}$, one can conclude that on its own, $C_{\{1,2,4,5,\cdots,p+1\}}$ forms a manifold structure of dimension $np-p(p-1)/2-(p-2)<np-p(p-1)/2$. Other cases $C_{\{i_1,\cdots,i_{p}\}}$ can also be considered the same way, and one still reaches the conclusion that, unless $\{i_1,\cdots,i_{p}\}=\{1,\cdots,p\}$, any of these $C_{\{i_1,\cdots,i_{p}\}}$ forms a manifold of dimension strictly less than $np-p(p-1)/2$.

\subsubsection{Density} \label{subsubsection 9.3.3}

\noindent
The discussion on dimensions above gives us a hint that as a subset, $S(n,p)^{*}=C_{\{1,\cdots,p\}}$ should be ``dense" on $S(n,p)$, as both can be endowed manifold structures of the same dimension $np-p(p-1)/2$ while other subsets $C_{\{i_1,\cdots,i_{p}\}}$ can only be manifolds of strictly lower dimensions. Let $M\in S(n,p)$ and $N$ be its reduced Cholesky factor. Then $N\in\mathbb{R}^{n\times p}_{*}$ is a mock lower triangular matrix whose mock diagonal entries $N_{ii}, i=1,\cdots,p$ are either all positive, in the case $M\in S(n,p)^{*}$, or zeros in some places, otherwise. Suppose, for concreteness, $M\in C_{\{1,2,4,5,\cdots,p+1\}}$ and so $N$ has the form as in \ref{ex2}. We can easily approximate (in Frobenius sense) $M$ with $M_{\epsilon}\in S(n,p)^{*}$ where $M_{\epsilon}=N_{\epsilon}N_{\epsilon}^{T}$ and,
\begin{equation} \label{ex3} 
    N_{\epsilon}=\begin{pmatrix}a_{11}&0&0&\dots&0&0\\ a_{21}&a_{22}&0&\dots&0&0\\a_{31}&a_{32}&\epsilon&\dots&0&0\\a_{41}&a_{42}&a_{43}&\dots&0&0\\
    \vdots&\vdots&\vdots&\ddots&\vdots&\vdots\\ a_{p1}&a_{p2}&a_{p3}&\dots&a_{p\,p-1}&\epsilon\\ a_{p+1\,1}&a_{p+1\,2}&a_{p+1\,3}&\dots&a_{p+1\,p-1}&a_{p+1\,p}\\ \vdots&\vdots&\vdots&\dots&\vdots&\vdots\\ a_{n1}&a_{n2}&a_{n3}&\dots&a_{n\,p-1}&a_{np}\end{pmatrix},
\end{equation}
for some $\epsilon>0$. In other words, $N_{\epsilon}$ is obtained from $N$ by filling in zero mock diagonal entries with $\epsilon$. Take $n=4,p=3$, then $N,N_{\epsilon}$ in \ref{ex2},\ref{ex3}, respectively, become:
\begin{equation*}
    N=\begin{pmatrix}a&0&0\\b&e&0\\c&f&0\\d&g&h\end{pmatrix}\,\,\text{ and }\,\,N_{\epsilon}=\begin{pmatrix}a&0&0\\b&e&0\\c&f&\epsilon\\d&g&h\end{pmatrix}.
\end{equation*}
Their corresponding $S(n,p)$ matrices are:
\begin{equation*}
    M =\begin{pmatrix}a^2&ab&ac&ad\\ab&b^2+e^2&bc+ef&bd+eg\\ac&bc+ef&c^2+f^2&cd+fg\\ad&bd+eg&cd+fg&d^2+g^2+h^2\end{pmatrix}\,\, \text{ and } \,\,M_{\epsilon}=\begin{pmatrix}a^2&ab&ac&ad\\ab&b^2+e^2&bc+ef&bd+eg\\ac&bc+ef&c^2+f^2+\epsilon^2&cd+fg+\epsilon h\\ad&bd+eg&cd+fg+\epsilon h&d^2+g^2+h^2\end{pmatrix}.
\end{equation*}
It's clear that $\|M-M_{\epsilon}\|_{F}\lesssim\|M\|_{\infty}\epsilon$ or $\|M-M_{\epsilon}\|_{F}= O_{M}(\epsilon)$. \\
(This also explains the choice of $M_{\epsilon} = \begin{pmatrix}1&1&1\\ 1&1+\epsilon^2&1+\epsilon\\ 1&1+\epsilon&2\end{pmatrix}$ for $M=\begin{pmatrix}1&1&1\\1&1&1\\1&1&2\end{pmatrix}$ in section \ref{section 2}.) \\
The point of this discussion is as follows. Consider $S(n,p)\subset\mathbb{R}^{n\times n}$ with the usual Frobenius metric. Let $M\in S(n,p)-S(n,p)^{*}$. Then there is $\epsilon>0$ small enough, such that within $B_2(M,\epsilon)\cap S(n,p)$, there exists $M^{*}\in S(n,p)^{*}\bigcap\big(B_2(M,\epsilon)\cap S(n,p)\big)$.

\subsection{A fact} \label{subsection 9.4}

\noindent Let $N\in\mathcal{L}^{n\times p}_{*}$ and $X\in\mathcal{L}^{n\times p}$. If $NX^{T}$ is a skew symmetric matrix, then $X=0_{n\times p}$.

\begin{proof}
That $NX^{T}$ is a skew symmetric matrix means that $(NX^{T})_{ii}=0, i=1,\cdots,n$. Let $Y= X^{T}$. Since $X_{ij}=0$ if $i<j$, $Y_{ij}=0$ if $i>j$. Now 
\begin{equation*} 
    (NX^{T})_{11} = (NY)_{11} = N_{11}Y_{11}=0\Rightarrow Y_{11}=0
\end{equation*} 
as $N_{11}>0$. That sets $col(Y)_1=\vec{0}$ and hence $col(NY)_1 =\vec{0}$ as well $row(NY)_1=\vec{0}$. This last fact sets $row(Y)_1=\vec{0}$, as now
\begin{equation*}
    0=(NY)_{1j}=\langle col(Y)_{j},row(N)_1\rangle = N_{11}Y_{1j}\Rightarrow Y_{1j}=0
\end{equation*}
as again, $N_{11}>0$, and $j=1,\cdots,n$. This now means 
\begin{equation*} 
    0=(NY)_{22}=N_{22}Y_{22}\Rightarrow Y_{22}=0
\end{equation*}
as $N_{22}>0$. Hence $col(Y)_2=\vec{0}$ and $col(NY)_2=\vec{0}$, $row(NY)_2=\vec{0}$.\\

\noindent
The rest of the proof proceeds as an induction scheme: 
\begin{center}
    ``If $col(Y)_{j}=\vec{0}$ $\big($then $col(NY)_{j}=\vec{0}$ and $row(NY)_{j}=\vec{0}\big)$ then $col(Y)_{j+1}=\vec{0}$."
\end{center} 
The first part of the proof is the base of this induction. The induction steps follow similarly to the base case. If $row(NY)_{j}=\vec{0}$ then $0=(NY)_{j\,j+1}=N_{jj}Y_{j\,j+1}$, as other terms are already determined to be zeros in the previous steps. Hence $Y_{j\,j+1}=0$ as $N_{jj}>0$. Since $0=(NY)_{j+1\,j+1}=N_{j+1\,j+1}Y_{j+1\,j+1}$, again, for the same reason, $Y_{j+1\,j+1}=0$. Hence $col(Y)_{j+1}=\vec{0}$.
\end{proof}

\subsection{Computational complexity of algorithms considered in section \ref{section 7}}

\noindent
Given $m$ different $n \times n$ matrices $\{ B_{i} \}_{i = 1}^m$. These are ``noisy", usually full-rank, versions of some low-rank matrix (assumed rank $p\ll n$) that needs to be unearthed through some averaging process. All three algorithms, LRC, dPCA, EigenV-ave, start with performing eigendecompositions for
\begin{equation} \label{eigdec}
    B_{i} = U_{i} \Lambda_{i} U_{i}^T \quad i=1,\cdots, m,
\end{equation}

\noindent
{\bf LRC}. The $\Lambda_{i}^{1/2} U_{i}$'s obtained from \ref{eigdec} are used to form the Cholesky factors $L_{i}$'s via \ref{QR}. Their Cholesky mean $\mathcal{F}(L_1,\cdots,L_m)$ is then used to give $\bar{B}_{LRC} = \mathcal{F}(L_1,\cdots,L_m)\mathcal{F}(L_1,\cdots,L_m)^{T}$. The mean column space is obtained by taking the top $p$ eigenvectors of $\bar{B}_{LRC}$.\\

\noindent
{\bf dPCA.} The $U_{i}$'s obtained from \ref{eigdec} are used to get $\tilde{B}_{i} = U_{i}U_{i}^T\, i=1,\cdots, m,$, from which their mean is calculated: $\bar{B}_{dPCA} = \frac{\sum_{i = 1}^m \tilde{B}_{i}}{m}$. The mean column space is obtained by taking the top $p$ eigenvectors of $\bar{B}_{dPCA}$.\\

\noindent
{\bf EigenV-ave.} The $U_{i}$'s obtained from \ref{eigdec} are averaged: $\bar{U} = \frac{\sum_{i = 1}^m  U_{i}}{m}$. The mean column space is obtained by the top $p$ eigenvectors of $\bar{U}\bar{U}^{T}$.\\

\noindent
In summary, all three use $(m+1)$ eigendecompositions on $n\times n$ matrices ($O((m+1)n^3)$). In addition, LRC performs $m$ QR decompositions on $p\times n$ matrices ($ O(m(p^3+p^3(n-p)))$), one summation of $m$ $n\times p$ matrices ($O(mnp)$), and one matrix multiplication of $n\times p$ matrices ($O(n^2p)$), while dPCA performs $m$ matrix multiplications of $n\times p$ matrices ($O(mn^2p)$) and one summation of $m$ $n\times n$ matrices ($O(mn^2)$), and lastly, EigenV-ave performs one summation of $m$ $n\times p$ matrices ($O(mnp)$). Note that the computational cost of the QR step \ref{QR} is $O(p^3+p^3(n-p))$, as it involves solving a QR decomposition of a $p\times p$ matrix and $n-p$ consistent linear systems of size $p\times p$.

\noindent
To facilitate a more concrete comparison, we give the empirical computation costs for the three algorithms, using $\Sigma$-perturbed random matrices, in the following setting. We first create $W\in\mathbb{R}^{n \times 3}$ as in section \ref{section 7} and let
\begin{equation*}
    \Sigma = W\Lambda_3 W^T,
\end{equation*}
where $\Lambda_3 = diag(10, 5, 2.5)$. We then generate $m$ noisy symmetric matrices
\begin{equation*}
    \Sigma^{(i)} = \Sigma + [E^{(i)} + (E^{(i)})^{T} ] /2 \quad i=1,\cdots, m,
\end{equation*}
where $E^{(i)}$ is a random matrix such that $E^{(i)}_{(l,j)} \sim N(0, 0.3^2)$. We consider $n \in \{100, 200, 500 \}$ and $m = 100$. The following results are generated from a personal laptop with 16GB RAM running the MacOS Monterey Version 12.1 system. Table \ref{tab:CompTime} based on 100 replicates shows that the computation time for all algorithms increases with dimension $n$. The computation time for LRC is slightly longer than EigenV-ave but is shorter than dPCA because LRC avoids multiple matrix multiplications before taking average.

\begin{table}[H]	
	{	
\textcolor{blue}{		\setlength{\tabcolsep}{5pt}
		\def\arraystretch{1.25}
		\begin{center}
			\begin{small}
				\begin{tabular}{cccc}
					Dimension (n) & LRC & dPCA & EigenV-ave \\[.1cm]
				\hline
					100 & 0.37 (0.07) & 0.44 (0.07) & 0.34 (0.06) \\
					200 & 1.81 (0.21)  & 2.22 (0.20) & 1.76 (0.23)  \\ 
					500 &  19.70 (2.01)& 22.19 (1.31) & 19.45 (2.04)   \\ 
				\end{tabular}
			\end{small}
		\end{center}
	}}
	\caption{Summary of mean (standard deviation) computational time in seconds on $\Sigma$-perturbed random matrices}
	\label{tab:CompTime}	
\end{table}


\begin{thebibliography}{10}

\bibitem{absil2009geometric}
P-A Absil, Mariya Ishteva, Lieven De~Lathauwer, and Sabine Van~Huffel.
\newblock A geometric newton method for oja's vector field.
\newblock {\em Neural computation}, 21(5):1415--1433, 2009.

\bibitem{absil2009optimization}
P-A Absil, Robert Mahony, and Rodolphe Sepulchre.
\newblock {\em Optimization algorithms on matrix manifolds}.
\newblock Princeton University Press, 2009.

\bibitem{adler2002newton}
Roy~L Adler, Jean-Pierre Dedieu, Joseph~Y Margulies, Marco Martens, and Mike
  Shub.
\newblock Newton's method on riemannian manifolds and a geometric model for the
  human spine.
\newblock {\em IMA Journal of Numerical Analysis}, 22(3):359--390, 2002.

\bibitem{ando2004geometric}
Tsuyoshi Ando, Chi-Kwong Li, and Roy Mathias.
\newblock Geometric means.
\newblock {\em Linear algebra and its applications}, 385:305--334, 2004.

\bibitem{bakker2018dynamic}
Craig Bakker, Mahantesh Halappanavar, and Arun~Visweswara Sathanur.
\newblock Dynamic graphs, community detection, and riemannian geometry.
\newblock {\em Applied network science}, 3(1):1--30, 2018.

\bibitem{bonnabel2013rank}
Silvere Bonnabel, Anne Collard, and Rodolphe Sepulchre.
\newblock Rank-preserving geometric means of positive semi-definite matrices.
\newblock {\em Linear Algebra and its Applications}, 438(8):3202--3216, 2013.

\bibitem{bonnabel2010riemannian}
Silvere Bonnabel and Rodolphe Sepulchre.
\newblock Riemannian metric and geometric mean for positive semidefinite
  matrices of fixed rank.
\newblock {\em SIAM Journal on Matrix Analysis and Applications},
  31(3):1055--1070, 2010.

\bibitem{bryner2017endpoint}
Darshan Bryner.
\newblock Endpoint geodesics on the stiefel manifold embedded in euclidean
  space.
\newblock {\em SIAM Journal on Matrix Analysis and Applications},
  38(4):1139--1159, 2017.

\bibitem{davidson1996c}
Kenneth~R Davidson.
\newblock {\em C*-algebras by example}, volume~6.
\newblock American Mathematical Soc., 1996.

\bibitem{dummit2004abstract}
David~Steven Dummit and Richard~M Foote.
\newblock {\em Abstract algebra}, volume~3.
\newblock Wiley Hoboken, 2004.

\bibitem{fan2019distributed}
Jianqing Fan, Dong Wang, Kaizheng Wang, and Ziwei Zhu.
\newblock Distributed estimation of principal eigenspaces.
\newblock {\em Annals of statistics}, 47(6):3009, 2019.

\bibitem{faraki2016image}
Masoud Faraki, Mehrtash~T Harandi, and Fatih Porikli.
\newblock Image set classification by symmetric positive semi-definite
  matrices.
\newblock In {\em 2016 IEEE Winter conference on applications of computer
  vision (WACV)}, pages 1--8. IEEE, 2016.

\bibitem{faraut1994analysis}
Jacques Faraut.
\newblock Analysis on symmetric cones.
\newblock {\em Oxford mathematical monographs}, 1994.

\bibitem{gamelin2003complex}
Theodore Gamelin.
\newblock {\em Complex analysis}.
\newblock Springer Science \& Business Media, 2003.

\bibitem{gentle2012numerical}
James~E Gentle.
\newblock {\em Numerical linear algebra for applications in statistics}.
\newblock Springer Science \& Business Media, 2012.

\bibitem{golub2013matrix}
Gene~H Golub and Charles~F Van~Loan.
\newblock {\em Matrix computations}, volume~3.
\newblock JHU press, 2013.

\bibitem{guillemin2010differential}
Victor Guillemin and Alan Pollack.
\newblock {\em Differential topology}, volume 370.
\newblock American Mathematical Soc., 2010.

\bibitem{Guo:2010}
Kai Guo, Prakash Ishwar, and Janusz Konrad.
\newblock Action recognition in video by sparse representation on covariance
  manifolds of silhouette tunnels.
\newblock In {\em Proceedings of the 20th International Conference on
  Recognizing Patterns in Signals, Speech, Images, and Videos}, pages 294--305,
  2010.

\bibitem{helgason1979differential}
Sigurdur Helgason.
\newblock {\em Differential geometry, Lie groups, and symmetric spaces}.
\newblock Academic press, 1979.

\bibitem{helmke1995critical}
Uwe Helmke and Mark~A Shayman.
\newblock Critical points of matrix least squares distance functions.
\newblock {\em Linear Algebra and its Applications}, 215:1--19, 1995.

\bibitem{hosseini2015matrix}
Reshad Hosseini and Suvrit Sra.
\newblock Matrix manifold optimization for gaussian mixtures.
\newblock {\em Advances in Neural Information Processing Systems}, 28:910--918,
  2015.

\bibitem{huber2004robust}
Peter~J Huber.
\newblock {\em Robust statistics}, volume 523.
\newblock John Wiley \& Sons, 2004.

\bibitem{journee2010low}
Michel Journ{\'e}e, Francis Bach, P-A Absil, and Rodolphe Sepulchre.
\newblock Low-rank optimization on the cone of positive semidefinite matrices.
\newblock {\em SIAM Journal on Optimization}, 20(5):2327--2351, 2010.

\bibitem{koch2007dynamical}
Othmar Koch and Christian Lubich.
\newblock Dynamical low-rank approximation.
\newblock {\em SIAM Journal on Matrix Analysis and Applications},
  29(2):434--454, 2007.

\bibitem{kolberg2020co}
Liis Kolberg, Nurlan Kerimov, Hedi Peterson, and Kaur Alasoo.
\newblock Co-expression analysis reveals interpretable gene modules controlled
  by trans-acting genetic variants.
\newblock {\em Elife}, 9:e58705, 2020.

\bibitem{lanckriet2004learning}
Gert~RG Lanckriet, Nello Cristianini, Peter Bartlett, Laurent~El Ghaoui, and
  Michael~I Jordan.
\newblock Learning the kernel matrix with semidefinite programming.
\newblock {\em Journal of Machine learning research}, 5(Jan):27--72, 2004.

\bibitem{lee2009manifolds}
Jeffrey Lee and Jeffrey~Marc Lee.
\newblock {\em Manifolds and Differential Geometry}, volume 107.
\newblock American Mathematical Soc., 2009.

\bibitem{li2019modeling}
Lingge Li, Dustin Pluta, Babak Shahbaba, Norbert Fortin, Hernando Ombao, and
  Pierre Baldi.
\newblock Modeling dynamic functional connectivity with latent factor gaussian
  processes.
\newblock {\em Advances in neural information processing systems},
  32:8263--8273, 2019.

\bibitem{li2014conformational}
Xiao-Bo Li and Forbes~J Burkowski.
\newblock Conformational transitions and principal geodesic analysis on the
  positive semidefinite matrix manifold.
\newblock In {\em International Symposium on Bioinformatics Research and
  Applications}, pages 334--345. Springer, 2014.

\bibitem{liang2014improved}
Yingyu Liang, Maria-Florina Balcan, Vandana Kanchanapally, and David~P
  Woodruff.
\newblock Improved distributed principal component analysis.
\newblock In {\em NIPS}, 2014.

\bibitem{liesecke2018ranking}
Franziska Liesecke, Dimitri Daudu, Rodolphe~Dug{\'e} de~Bernonville,
  S{\'e}bastien Besseau, Marc Clastre, Vincent Courdavault, Johan-Owen
  De~Craene, Joel Cr{\`e}che, Nathalie Giglioli-Guivarc'h, Ga{\"e}lle
  Gl{\'e}varec, et~al.
\newblock Ranking genome-wide correlation measurements improves microarray and
  rna-seq based global and targeted co-expression networks.
\newblock {\em Scientific reports}, 8(1):1--16, 2018.

\bibitem{lin2019riemannian}
Zhenhua Lin.
\newblock Riemannian geometry of symmetric positive definite matrices via
  cholesky decomposition.
\newblock {\em SIAM Journal on Matrix Analysis and Applications},
  40(4):1353--1370, 2019.

\bibitem{massart2020quotient}
Estelle Massart and P-A Absil.
\newblock Quotient geometry with simple geodesics for the manifold of
  fixed-rank positive-semidefinite matrices.
\newblock {\em SIAM Journal on Matrix Analysis and Applications},
  41(1):171--198, 2020.

\bibitem{meyer2009subspace}
Gilles Meyer, Michel Journ{\'e}e, Silvere Bonnabel, and Rodolphe Sepulchre.
\newblock From subspace learning to distance learning: a geometrical
  optimization approach.
\newblock In {\em 2009 IEEE/SP 15th Workshop on Statistical Signal Processing},
  pages 385--388. IEEE, 2009.

\bibitem{mishra2011low}
Bamdev Mishra, Gilles Meyer, and Rodolphe Sepulchre.
\newblock Low-rank optimization for distance matrix completion.
\newblock In {\em 2011 50th IEEE Conference on Decision and Control and
  European Control Conference}, pages 4455--4460. IEEE, 2011.

\bibitem{orsi2006newton}
Robert Orsi, Uwe Helmke, and John~B Moore.
\newblock A newton-like method for solving rank constrained linear matrix
  inequalities.
\newblock {\em Automatica}, 42(11):1875--1882, 2006.

\bibitem{pennec2006intrinsic}
Xavier Pennec.
\newblock Intrinsic statistics on riemannian manifolds: Basic tools for
  geometric measurements.
\newblock {\em Journal of Mathematical Imaging and Vision}, 25(1):127--154,
  2006.

\bibitem{qu2002principal}
Yongming Qu, George Ostrouchov, Nagiza Samatova, and Al~Geist.
\newblock Principal component analysis for dimension reduction in massive
  distributed data sets.
\newblock In {\em Proceedings of IEEE International Conference on Data Mining
  (ICDM)}, volume 1318, page 1788, 2002.

\bibitem{schiratti2017bayesian}
Jean-Baptiste Schiratti, St{\'e}phanie Allassonni{\`e}re, Olivier Colliot, and
  Stanley Durrleman.
\newblock A bayesian mixed-effects model to learn trajectories of changes from
  repeated manifold-valued observations.
\newblock {\em The Journal of Machine Learning Research}, 18(1):4840--4872,
  2017.

\bibitem{shen2002hammer}
Dinggang Shen and Christos Davatzikos.
\newblock Hammer: hierarchical attribute matching mechanism for elastic
  registration.
\newblock {\em IEEE transactions on medical imaging}, 21(11):1421--1439, 2002.

\bibitem{shirer2012decoding}
William~R Shirer, Srikanth Ryali, Elena Rykhlevskaia, Vinod Menon, and
  Michael~D Greicius.
\newblock Decoding subject-driven cognitive states with whole-brain
  connectivity patterns.
\newblock {\em Cerebral cortex}, 22(1):158--165, 2012.

\bibitem{strichartz2000way}
Robert~S Strichartz.
\newblock {\em The way of analysis}.
\newblock Jones \& Bartlett Learning, 2000.

\bibitem{vandereycken2009embedded}
Bart Vandereycken, P-A Absil, and Stefan Vandewalle.
\newblock Embedded geometry of the set of symmetric positive semidefinite
  matrices of fixed rank.
\newblock In {\em 2009 IEEE/SP 15th Workshop on Statistical Signal Processing},
  pages 389--392. IEEE, 2009.

\bibitem{vandereycken2013riemannian}
Bart Vandereycken, P-A Absil, and Stefan Vandewalle.
\newblock A riemannian geometry with complete geodesics for the set of positive
  semidefinite matrices of fixed rank.
\newblock {\em IMA Journal of Numerical Analysis}, 33(2):481--514, 2013.

\bibitem{zeestraten2017approach}
Martijn~JA Zeestraten, Ioannis Havoutis, Joao Silv{\'e}rio, Sylvain Calinon,
  and Darwin~G Caldwell.
\newblock An approach for imitation learning on riemannian manifolds.
\newblock {\em IEEE Robotics and Automation Letters}, 2(3):1240--1247, 2017.

\end{thebibliography}
\end{document}